\documentclass[final,3p,times]{elsarticle}

\usepackage{amssymb}
\usepackage{amsmath}
 \usepackage{amsthm}

\usepackage[T1]{fontenc}
\usepackage{tikz, mdframed, subcaption, amsmath, multirow, tikz, graphicx, enumitem, lmodern}
\usetikzlibrary{tikzmark,calc}

\newtheorem{theorem}{Theorem}

\newtheorem{corollary}{Corollary}
\newtheorem{lemma}{Lemma}
\newtheorem{example}{Example}
\newtheorem{reduction}{Reduction}
\newtheorem{remark}{Remark}

\newcommand{\osd}{\gamma^{\rm{OD}}}
\newcommand{\od}{\gamma^{\rm{OD}}}
\newcommand{\otd}{\gamma^{\rm{OTD}}}
\newcommand{\old}{\gamma^{\rm{OTD}}}
\newcommand{\ld}{\gamma^{\rm{LD}}}

\newcommand{\ltd}{\gamma^{\rm{LTD}}}

\newcommand{\x}{\gamma^{\rm{X}}}

\newcommand{\OSD}{{\textsc{OD}}}
\newcommand{\OD}{{\textsc{OD}}}
\newcommand{\OTD}{{\textsc{OTD}}}
\newcommand{\OLD}{{\textsc{OTD}}}
\newcommand{\DTD}{{\textsc{ITD}}}
\newcommand{\ITD}{{\textsc{ITD}}}
\newcommand{\LTD}{{\textsc{LTD}}}
\newcommand{\ID}{{\textsc{ID}}}
\newcommand{\LD}{{\textsc{LD}}}
\newcommand{\X}{{\textsc{X}}}

\newcommand{\odadmis}{{\textsc{OD}}-admissible}
\newcommand{\otdadmis}{{\textsc{OTD}}-admissible}
\newcommand{\oldadmis}{{\textsc{OTD}}-admissible}

\newcommand{\hyp}{\mathcal{H}}
\newcommand{\clu}{\mathcal{C}}
\newcommand{\for}{\mathcal{F}^1}
\newcommand{\fac}{\mathcal{F}^2}
\newcommand{\F}{\mathcal{F}}

\newcommand{\NP}{{\textsf{NP}}}
\newcommand{\calC}{\mathcal{C}}
\newcommand{\calO}{\mathcal{O}}

\newcommand{\yes}{\textsc{Yes}}

\newcommand{\ddelta}{\bigtriangleup}
\newcommand{\mincode}[1]{\textsc{#1-Code}}

\newcommand{\defproblem}[3]{
\vspace{3mm}
\noindent \fbox{
\begin{minipage}{0.96\textwidth}
\begin{tabular*}{\textwidth}{@{\extracolsep{\fill}}lr}
#1 \\ \end{tabular*}
{\bf{Input:}} #2 \\
{\bf{Question:}} #3
\end{minipage}
}

\vspace{3mm}
}

\begin{document}

\begin{frontmatter}



\title{On open-separating dominating codes in graphs}


\author[label1,label2]{Dipayan Chakraborty} 
\author[label1]{Annegret K. Wagler} 

\affiliation[label1]{organization={Universit\'{e} Clermont-Auvergne, CNRS, Mines de Saint-\'{E}tienne, Clermont-Auvergne-INP, LIMOS},
            city={Clermont-Ferrand},
            postcode={63000}, 
            country={France}}
            
\affiliation[label2]{organization={Department of Mathematics and Applied Mathematics, University of Johannesburg},
            city={Auckland Park},
            postcode={2006}, 
            country={South Africa}}

\begin{abstract}
Using dominating sets to separate vertices of graphs is a well-studied problem in the larger domain of identification problems. In such problems, the objective is to choose a suitable dominating set $C$ of a graph $G$ which is also separating in the sense that the neighbourhoods of any two distinct vertices of $G$ have distinct intersections with $C$. Such a dominating and separating set $C$ of a graph is often referred to as a \emph{code} in the literature. Depending on the types of dominating and separating sets used, various problems arise under various names in the literature. In this paper, we introduce a new problem in the same realm of identification problems whereby the code, called  \emph{open-separating dominating code}, or \emph{\OD-code} for short, is a dominating set and uses open neighbourhoods for separating vertices. 
The paper studies the fundamental properties concerning the existence, hardness and minimality of \OSD-codes. Due to the emergence of a close and yet difficult to establish relation of the \OSD-code with another well-studied code in the literature called open (neighborhood)-locating dominating code (referred to as the \emph{open-separating total-dominating code} and abbreviated as \emph{\OLD-code} in this paper), we compare the two codes on various graph families. Finally, we also provide an equivalent reformulation of the problem of finding \OSD-codes of a graph as a covering problem in a suitable hypergraph and discuss the polyhedra associated with \OSD-codes, again in relation to \OLD-codes of some graph families already studied in this context.
\end{abstract}




\end{frontmatter}

\section{Introduction}
The problem of placing surveillance devices in buildings to locate an intruder (like a fire, a thief or a saboteur) leads naturally to different separation-domination type problems in the graph modeling the building (where rooms are represented as vertices and connections between rooms as edges). Depending on the characteristics of the detection devices (to detect an intruder only if it is present in the room where the detector is installed and/or to detect one in any neighbouring room), different kinds of dominating sets can be used to detect the existence of an intruder, whereas different separation-type properties are considered to exactly locate the position of an intruder in the building.
 
More precisely, let $G=(V(G), E(G)) = (V,E)$ be a graph and let $N(v) = N_G(v) = \{u \in V : uv \in E\}$ (respectively, $N[v] = N_G[v] = N(v) \cup \{v\}$) denote the \emph{open neighbourhood} (respectively, \emph{closed neighborhood}) of a vertex $v \in V$. A subset $C \subseteq V$ is \emph{dominating} (respectively, \emph{total-dominating}) if $N[v]\cap C$ (respectively, $N(v)\cap C$) is a non-empty set for each $v \in V$. Moreover, $C \subseteq V$ is \emph{closed-separating} (respectively, \emph{open-separating}) if $N[v]\cap C$ (respectively, $N(v)\cap C$) is a unique set for each $v \in V$. Furthermore, the set $C$ is \emph{locating} if $N(v)\cap C$ is a unique set for each $v \in V\setminus C$.

So far, the following combinations of separation and domination properties have been studied in the literature over the last decades:
\begin{itemize}[leftmargin=12pt, itemsep=0pt]
\item closed-separation with domination and total-domination leading to \emph{identifying codes} (\emph{\ID-codes} for short)~\cite{KCL1998} and 
\emph{identifying total-dominating codes}\footnote{Identifying total-dominating codes had been introduced to the literature in~\cite{HHH2006} under the name \emph{differentiating total-dominating codes}. However, due to consistency in notation, we prefer to call them ITD-codes in this article.} (\emph{\DTD-codes} for short)~\cite{HHH2006}, respectively;
\item location with domination and total-domination leading to \emph{locating dominating codes} (\emph{\LD-codes} for short)~\cite{Slater1988} and \emph{locating total-dominating codes} (\emph{\LTD-codes} for short)~\cite{HHH2006}, respectively;
\item open-separation with total-domination leading to 
\emph{open-separating total-do\-minating codes}\footnote{OTD-codes were introduced independently in~\cite{HLR2002} and in~\cite{SS2010} under the names of \emph{strongly $(t,\le l)$-identifying codes} and \emph{open neighborhood locating-dominating sets} (or \emph{OLD-sets}), respectively. However, due to consistency in naming that specifies the separation and the domination property, we prefer to call them open-separating total-dominating codes in this article.} (\emph{\OLD-codes} for short) \cite{HLR2002,SS2010}.
\end{itemize}
Such problems have several applications, for example, in fault-detection in multiprocessor networks~\cite{KCL1998}, logical definability of graphs~\cite{PVV06}, canonical labeling of graphs for the graph isomorphism problem~\cite{B80} and locating intruders in facilities using sensor networks~\cite{UTS04}, to name a few. 
Recall that to model the latter monitoring problems, a building is represented as a graph and the task is to place sensors at some vertices. 
Hereby, the subset of rooms equipped with sensors corresponds to a code in the graph, the domination properties guarantee that all rooms are monitored, whereas the separation properties ensure that the intruder can be exactly located. For monitoring systems in which each sensor can distinguish between the presence of an intruder at its own vertex and at a neighboring vertex, the separating property useful in pinpointing the intruder is location~\cite{Slater1988}. On the other hand, if some sensors lose the ability to distinguish between the presence of an intruder at their own vertex and a neighboring one, then the separating property required to locate the intruder is that of closed-separation~\cite{KCL1998}. Finally, if some sensor in the monitoring system is disabled by an intruder, then the separating property useful in locating the intruder is that of open-separation \cite{HLR2002,SS2010}. For more applications and results, we refer to an extensive internet bibliography containing over 500 articles around these topics maintained by Jean and Lobstein~\cite{JL_lib}.

\begin{figure}[h]
\begin{center}
\includegraphics[scale=1.0]{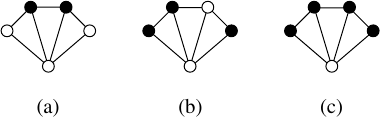}
\caption{Minimum X-codes in a graph (the black vertices belong to the code), where (a) 
is both an LD- and LTD-code, (b) an OD-code, (c) an OTD-, ID-, and ITD-code.}
\label{fig_exp_x-codes}
\end{center}
\end{figure}

In this paper, we aim at studying open-separation combined with domination. We call a subset $C \subseteq V$ of a graph $G=(V,E)$ an \emph{open-separating dominating code} (\emph{\OSD-code} for short) if it is
\begin{itemize}[leftmargin=12pt, itemsep=0pt]
\item a dominating set, that is $N[v]\cap C$ is a non-empty set for each $v \in V$; and
\item an open-separating set, that is $N(v)\cap C$ is a unique set for each $v \in V$.
\end{itemize}
As is evident from the definitions, the \OD-code is a variation of the \OTD-code which combines open-separation with total-domination. In the literature of domination-based identification problems, the \OTD-code has been quite well studied, for example, in~\cite{ArBiLuWa-20,ABLW_2022,CaCoEFo-22,ChFoPaWa-24,Che-14,FGRS2021,FoMeNaPaVa17a,FoMeNaPaVa17b,Giv-22,HeYe-14,JeSe-23,Kin-15,PP_2017,SS2010,SeSl11}. The \OD-code, however, has so far not been studied in the literature of identification problems except for in the conference proceeding~\cite{CW_ISCO2024} of the current paper where some of the same results have appeared without proofs. The motivation to combine open-separation with domination to study the \OD-codes is a most natural one, as has been the case with other separation properties like location and closed-separation which have also been studied in combination with both domination and total-domination. Figure \ref{fig_exp_x-codes} illustrates examples of such codes in a small graph.

Note that not all graphs admit codes of all the studied types. Accordingly, in Section~\ref{sec2} of this paper, we address the conditions for the existence of \OSD-codes and their relations to codes of other types. It turns out that \OSD-codes possess a particularly close relationship with \OLD-codes as the minimum cardinalities of the two differ by at most one. Moreover, for any $\X \in \{\ID, \DTD, \LD, \LTD, \OLD\}$, the problem of determining an \X-code of minimum cardinality $\gamma^{X}(G)$ of a graph $G$, called the \emph{\X-number} of $G$, has been shown to be \NP-hard~\cite{CHL2003,CS1987,SS2010}. 
In Section \ref{sec3}, we show that \NP-hardness holds for \OSD-codes as well. This also enables us to reprove the \NP-hardness of the \OTD-problem on a narrower class of input graphs and with an arguably simpler reduction than previously considered, for example, in~\cite{PP_2017,SS2010}. Furthermore, despite the close relationship between the \OSD- and \OLD-numbers of a graph, we show that deciding whether the two numbers differ is \NP-hard. This motivates us to compare the \OSD- and the \OLD-codes of graphs of different families in Section \ref{sec4} and to study their related polyhedra in Section \ref{sec5}. We close with some concluding remarks and lines of future research.

\subsection*{Further terminologies and notations}
A graph $H$ is called a \emph{subgraph} of a graph $G$ if $V(H) \subseteq V(G)$ and $E(H) \subseteq E(G)$. Moreover, a subgraph $H$ of $G$ is said to be \emph{induced} if $uv \in E(H)$ if and only if $uv \in E(G)$. For an induced subgraph $H$ of $G$, if $S = V(H)$, then $H$ is said to be \emph{induced by} $S$ and is denoted $G[S]$. The \emph{girth} of a given graph $G$ is the smallest length of an induced cycle in $G$. The \emph{degree} of a vertex $v$, denoted $\deg(v) = \deg_G(v)$, is the number of neighbors of $v$ in $G$. A vertex of degree~$0$ is called an \emph{isolated vertex}, a vertex of degree~$1$ is called a \emph{leaf} and the neighbor of a leaf is called the \emph{support vertex} of the leaf. 
Given any two vertices of $u,v$ of a graph $G$, the \emph{distance} between them, denoted $d(u,v) = d_G(u,v)$, is the length of the shortest induced path in $G$ with $u$ and $v$ as its endpoints. A graph is said to be \emph{complete} (on $n$ vertices), and denoted $K_n$, if any two of its vertices are adjacent. A \emph{clique} is a vertex subset $S$ of a graph $G$ such that the subgraph of $G$ induced by $S$ is complete. On the other hand, a \emph{stable set} is a vertex subset $S$ of a graph $G$ such that no two vertices in $S$ are adjacent in $G$. A graph $G=(U \cup W, E)$ is \emph{bipartite} if its vertex set can be partitioned into two stable sets $U$ and $W$ so that every edge of $G$ has one endpoint in $U$ and the other in $W$. A graph $G=(Q \cup S, E)$ is a \emph{split graph} if its vertex set can be partitioned into a clique $Q$ and a stable set $S$. Given two graphs $G_1$ and $G_2$ whose vertex sets are disjoint, the graph $(V(G_1) \cup V(G_2), E(G_1) \cup E(G_2))$ is denoted by $G_1 \oplus G_2$.

A vertex $w \in V$ is said to \emph{open-separate} two distinct vertices $u,v \in V$ if $w \in N(u) \ddelta N(v)$. Moreover, a vertex subset $S$ of $G$ is said to \emph{open-separate} the pair $u,v$ (by a vertex $w$) if $w \in S$ such that $w$ open-separates $u,v$. Then, it can be verified that a set $S$ has the open-separating property (defined before) if and only if it open-separates every pair of distinct vertices of $G$.

\section{Reformulation, existence and bounds}
\label{sec2}

In this section, we reformulate the \OD-problem, address fundamental questions concerning the existence of \OSD-codes and bounds for \OSD-numbers. 

\subsection{Hypergraph representation of the \OSD-problem} 

Given a graph $G=(V,E)$ and a problem \X, we look for a hypergraph $\hyp_\X(G) = (V,\F_\X)$ so that $C \subset V$ is an \X-code of $G$ if and only if $C$ is a \emph{cover} of $\hyp_\X(G)$ satisfying $C \cap F \neq \emptyset$ for all $F \in \F_\X$. 
Then the \emph{covering number} $\tau(\hyp_\X(G))$, defined as the minimum cardinality of a cover of $\hyp_\X(G)$, equals by construction the \X-number $\x(G)$ of $G$. The hypergraph $\hyp_\X(G)$ is called the \emph{\X-hypergraph} of the graph $G$.

It is a simple observation that for an \X-problem involving domination (respectively, total-domination), $\F_\X$ needs to contain the closed (respectively, open) neighborhoods of all vertices of $G$. 
In order to encode the separation of vertices, that is, the fact that the intersections of an \X-code $C$ with the neighborhood of each vertex is \emph{unique}, it was suggested in \cite{ABLW_2018,ABLW_2022} to use the symmetric differences of the neighborhoods. 
Here, given two sets $A$ and $B$, their \emph{symmetric difference} is defined by $A \Delta B = (A \setminus B) \cup (B \setminus A)$. 
In fact, it has been shown in \cite{ABLW_2018,ABLW_2022} that a code $C$ of a graph $G$ is closed-separating (respectively, open-separating) if and only if $(N[u] \Delta N[v]) \cap C \neq \emptyset$ (respectively, $(N(u) \Delta N(v)) \cap C \neq \emptyset$) for all pairs of distinct vertices $u,v$ of $G$. This implies for \OSD-codes:

\begin{corollary}\label{cor_OSD-hyper}
The \OSD-hypergraph $\hyp_{\OSD}(G) = (V,\F_{\OSD})$ of a graph $G=(V,E)$ is composed of
\begin{itemize}[leftmargin=12pt, itemsep=0pt]
\item the closed neighborhoods $N[v]$ of all vertices $v \in V$ and 
\item the symmetric differences $N(u) \Delta N(v)$ of open neighborhoods of distinct vertices $u,v \in V$
\end{itemize}
as hyperedges in $\F_{\OSD}$ and $\osd(G)=\tau(\hyp_{\OSD}(G))$ holds.
\end{corollary}

\subsection{Existence of \OSD-codes}
It has been observed in the literature that the studied domination and separation properties may not apply to all graphs, see for example \cite{KCL1998,HHH2006,SS2010}. More precisely, 
\begin{itemize}[leftmargin=12pt, itemsep=0pt]
\item total-domination excludes the occurrence of \emph{isolated vertices} in graphs, that is, vertices $v$ with $N(v) = \emptyset$;
\item closed-separation (respectively, open-separation) excludes the occurrence of \emph{closed twins} (respectively, \emph{open twins}), that is, distinct vertices $u,v$ with $N[u] = N[v]$ (respectively, $N(u) = N(v)$).
\end{itemize}
Calling a graph $G$ to be \emph{\X-admissible} if $G$ has an \X-code, we see that while, for example, every graph $G$ is \LD-admissible, a graph $G$ is \OLD-admissible if and only if $G$ has neither isolated vertices nor open twins. Accordingly, we conclude the following regarding the existence of \OSD-codes of graphs.
\begin{corollary}\label{cor_OSD-admissible}
A graph $G$ is \OSD-admissible if and only if $G$ has no open twins. 
\end{corollary}
Since any two distinct isolated vertices of a graph are open twins with the empty set as both their open neighbourhoods, Corollary~\ref{cor_OSD-admissible} further implies that an \OSD-admissible graph has at most one isolated vertex.

\subsection{Bounds on \OSD-numbers and their relations to other \X-numbers.}

\begin{remark} \label{lem_unique open signature}
Let $G$ be an \OD-admissible graph and let $C$ be an \OD-code of $G$. Then, there exists at most one vertex $w$ of $G$ such that $N(w) \cap C=\emptyset$.
\end{remark}

\begin{proof}
Toward contradiction, if there exist two distinct vertices $u$ and $v$ of $G$ such that $N(u) \cap C = N(v) \cap C = \emptyset$, then $C$ does not open-separate the pair $u,v$, a contradiction. This proves the result.
\end{proof}

\begin{remark} \label{lem_char_osd}
Let $G$ be an \OD-admissible graph and let $C$ be a dominating set of $G$ such that there exists at most one vertex $w$ of $G$ with $N(w) \cap C=\emptyset$. If the set $C$ open-separates every pair $u,v$ of distinct vertices of $G$ with $d(u,v) \leq 2$, then $C$ is an \OD-code of $G$.
\end{remark}

\begin{proof}
It is enough to show that $C$ open-separates every pair $u,v$ of distinct vertices of $G$ such that $d(u,v) \geq 3$. So, assume that $u$ and $v$ are two such vertices of $G$ with $d(u,v) \geq 3$. Then, we have $N(u) \triangle N(v) = N(u) \cup N(v)$. If $C$ does not open-separate the pair $u,v$, then it implies that $C$ does not intersect $N(u) \cup N(v)$, that is, $N(u) \cap C = N(v) \cap C = \emptyset$, which contradicts our assumption. Therefore, $C$ open-separates $u,v$ and this proves the result.
\end{proof}

We now prove general upper and lower bounds on the \OSD-number of a graph in terms of its order, that is the number of vertices. The proof for these bounds involves the following result due to Bondy.

\begin{theorem}[Bondy~\cite{Bondy1972}]\label{thm:bondy}
Let $C$ be a set with $|C| = n$ and let $C_1, C_2, \ldots , C_n$ be subsets of $C$ for each $i \in [n]$ such that $C_i \neq C_j$ for each $i,j \in [n]$ and $i \neq j$. Then, there exists an element $x$ of $S$ such that $C_i \setminus \{x\} \neq C_j \setminus \{x\}$ for each $i,j \in [n]$ and $i \neq j$. 
\end{theorem} 

\begin{lemma} \label{lem:X-code S1+V2}
Let $G_1$ and $G_2$ be \OD-admissible graphs on disjoint vertex sets each of whose cardinalities is at least two. Let $G = G_1 \oplus G_2$ be their disjoint union. If $S_1$ is an \OD-code of $G_1$, then the set $S = S_1 \cup V(G_2)$ is an \OD-code of $G$.
\end{lemma}

\begin{proof}
Clearly, the set $S$ is a dominating set of $G$ since $S_1$ is a dominating set of $G_1$. Let $u,v$ be two distinct vertices of $G$. Then, we show that the set $S$ has non-empty intersection with $N_G(u) \ddelta N_G(v)$. If $u,v \in V(G_1)$, then we have $N_G(u) \ddelta N_G(v) = N_{G_1}(u) \ddelta N_{G_1}(v)$. Since $(N_{G_1}(u) \ddelta N_{G_1}(v)) \cap S_1 \ne \emptyset$ due to $S_1$ being an \OD-code of $G_1$, we have $(N_G(u) \ddelta N_G(v)) \cap S \ne \emptyset$. Now if $u,v \in V(G_2)$, then $N_G(u) \ddelta N_G(v) = N_{G_2}(u) \ddelta N_{G_2}(v) \ne \emptyset$ since $G_2$ is \X-admissible. This implies that $(N_{G_2}(u) \ddelta N_{G_2}(v)) \cap V(G_2) \ne \emptyset$ and hence, $(N_G(u) \ddelta N_G(v)) \cap S \ne \emptyset$. Finally, let us assume that $u \in V(G_1)$ and $v \in V(G_2)$. Since $G_2$ is on at least two vertices, this implies that there exists a vertex $v' \in N_{G_2}(v) \subseteq S$ with $v' \ne v$. Now, since $N_G(u) \ddelta N_G(v) = N_{G_1}(u) \cup N_{G_2}(v)$, we have $v' \in (N_G(u) \cup N_G(v)) \cap S$. This proves the result.
\end{proof}

\begin{theorem} \label{thm_osd_bounds-n}
For a graph $G$ on $n \geq 2$ vertices without open twins and any isolated vertices, we have $\lceil \log n \rceil  \leq \osd(G) \leq n-1$.
\end{theorem}

\begin{proof}
Let us first prove the result on the lower bound of $\od(G)$. Let $S$ be a minimum OD-code of $G$, that is, $|S| = \od(G)$. Then, for each vertex $v \in V(G)$, the subset $N(v) \cap S \subset S$ is unique. In other words, the number $n$ of vertices of $G$ can be at most the number of subsets of $S$ which is $2^{|S|}$. This implies that $\log n \le |S| = \od(G)$ which further implies $\lceil \log n \rceil \le \od(G)$, since the \OD-number is an integer.

We now prove the result on the upper bound of $\od(G)$. Since $G$ is without isolated vertices, it must have at least one edge. So, let $G_1$ be a component of $G$ with at least one edge and of order $n_1$. Then, $G_1$ is also \OD-admissible. Let $G_2$ be the disjoint union of all components of $G$ other than $G_1$. In other words, we have $G = G_1 \oplus G_2$. Let $S_1$ be a minimum \OD-code of $G_1$, that is, $|S_1| = \od(G_1)$. Now, if $\od(G_1) \le n_1-1$, then, by applying Lemma~\ref{lem:X-code S1+V2}, the set $S = S_1 \cup V(G_2)$ is an \OD-code of $G$. This implies that $\od(G_1) \le |S| \le n-1$.

Therefore, let us assume that $\od(G_1) = |S_1| = n_1$, that is, $S_1 = V(G_1)$. Now, let $V(G_1)=\{v_1, v_2, \ldots v_{n_1}\}$ and let $C_i = N(v_i)$ for each $i \in [n_1]$. Since $S_1$ is an \OD-code of $G_1$, we have $C_i \neq C_j$ for all $i, j \in [n_1]$ and $i \neq j$. However, by Theorem~\ref{thm:bondy} (Bondy's theorem), there exists a vertex $x_1$ of $G_1$ such that $C_i \setminus \{x_1\} \neq C_j \setminus \{x_1\}$ for all $i, j \in [n_1]$ and $i \neq j$. In other words, the set $S_1 \setminus \{x_1\}$ is an open-separating set of $G_1$. Moreover, $S_1\setminus \{x_1\}$ is also a dominating set of $G_1$, since $S_1 \setminus \{x_1\} = V(G_1) \setminus \{x_1\}$ and $G_1$ is connected. Thus, $S_1 \setminus \{x_1\}$ is an \OD-code of $G_1$ which is a contradiction to the minimality of $S_1$ in being an \OD-code of $G_1$. Therefore, we have $\od(G_1) \leq n_1-1$ and we are again done by Lemma~\ref{lem:X-code S1+V2}.
\end{proof}

Both bounds in Theorem~\ref{thm_osd_bounds-n} are tight. For the lower bound, an example of a graph whose \OD-number attains the bound would be a split graph, say $G=(K \cup S, E)$, where $K \cong K_{2n+1}$ and $S$ is a stable set such that for all subset $U$ of $V(K)$ with $|U| \notin \{0,2n\}$, there exists a vertex $v \in S$ with $N(v) = U$. Examples of graphs whose \OD-number attains the upper bound in Theorem~\ref{thm_osd_bounds-n} are addressed in Section~\ref{sec4}.

Next, we turn to the relation of \OD-numbers to other \X-numbers for $\X \in \{\LD, \LTD, \ID, \ITD, \OTD\}$. 
The next theorem brings out a close association between the \OD- and the \OTD-numbers of an \OTD-admissible graph.

\begin{theorem} \label{thm_osd_bounds-X}
Let $G$ be an \OD-admissible graph. 
\begin{itemize}[leftmargin=18pt, itemsep=0pt]
\item[(a)] If $G$ is a disjoint union of a graph $G'$ and an isolated vertex, then we have $\od (G) = \otd (G') + 1$.
\item[(b)] If $G$ is \OTD-admissible, we have $\otd (G) - 1 \leq \od (G) \leq \otd(G)$.
\end{itemize}
\end{theorem}

\begin{proof}
Let us first assume that $G = G' \oplus \{v\}$ for some vertex $v$ of $G$. Moreover, let $C$ be any minimum \OD-code of $G$. Then, $v \in C$ in order for $C$ to dominate $v$. This implies that $N(v) \cap C = \emptyset$. Therefore, by Remark~\ref{lem_unique open signature}, the set $C \setminus \{v\}$ must total-dominate all vertices in $V \setminus \{v\}$, that is, $C$ is an \OD-code of $G'$. This implies that $\otd(G') \leq |C| - 1 = \od(G)-1$, that is, $\od(G) \leq \otd(G') + 1$. On the other hand, if $C'$ is an \OTD-code of $G'$, then $C' \cup \{v\}$ is an \OD-code of $G$. This implies that $\od(G) \leq |C'|+1 = \otd(G')+1$ and hence, assertion (a) follows.

Let us now assume that $G$ has no isolated vertices. Then, $G$ is also \OTD-admissible. Now, the inequality $\od (G) \leq \otd (G)$ holds by the fact that any \OTD-code of $G$ is also a dominating and an open-separating set of $G$ and hence, is an \OD-code of $G$. To prove the other inequality, we let $C$ be a minimum \OD-code ~of $G$. Then, we have $|C| = \od(G)$. Now, if $C$ is also an \OTD-code of $G$, then the result holds trivially. Therefore, let us assume that $C$ is not an \OTD-code of $G$, that is, in particular, $C$ is not a total-dominating set of $G$. This implies that there exists a vertex $v$ of $G$ such that $C \cap N(v) = \emptyset$. Moreover, by Remark~\ref{lem_unique open signature}, $C \cap N(v)=\emptyset$ for exactly the vertex $v$ of $G$. This implies that $C \cup \{v\}$ is a total-dominating set of $G$ and hence, is also an \OTD-code of $G$. Therefore, we have $\otd (G) \leq |C|+1 = \od (G) +1$ which implies assertion (b).
\end{proof}

We next address the relation of \OD-numbers to \LD- and \LTD-numbers of a graph.
\begin{theorem}\label{thm_osd_bounds-L}
Let $G$ be an \OD-admissible graph. 
\begin{itemize}[leftmargin=18pt, itemsep=0pt]
\item[(a)] We have $\ld (G) \le \od (G)$.
\item[(b)] If $G$ has no isolated vertex, then we have $\ltd (G) - 1 \leq \od (G)$.
\end{itemize}
\end{theorem}

\begin{proof}
We clearly have, by definition, that every open-separating set of a graph $G$ is also a locating set of $G$. 
Hence, every \OD-code of $G$ is also an \LD-code of $G$ and assertion (a) follows. 
Moreover, if $G$ has no isolated vertex, then every \OTD-code of $G$ is also an \LTD-code of $G$ which implies that $\ltd (G) \leq \otd (G)$ holds.
Combining the latter inequality with Theorem~\ref{thm_osd_bounds-X}(b), we further conclude
$$ \ltd (G) - 1 \leq \otd (G)- 1 \leq \od (G) $$
which implies assertion (b).
\end{proof}

On the other hand, there exist graphs $G$ with $\od (G) < \ltd (G)$ and the \OSD-numbers of graphs are not generally comparable to their \ID- and \ITD-numbers, see Table~\ref{tab_comparisons} for an illustration. 

\begin{table}[!t]
\centering
\begin{tabular}{|c|c|c|c|c|c|}
\hline
\bf \X & \bf \X-code & \bf \OSD-code & \bf ~Graph name~ & ~$\mathbf{\osd(G)}$~ & ~$\mathbf{\x(G)}$~ \\ \hline \hline
\multirow{2}{*}{\X=\ID} & \begin{tikzpicture}[myv/.style={circle, draw, inner sep=1.3pt}, scale=1]
    \node (o) at (0,0) {};
    \node[myv, fill=black] (v1) at (0,-0.25) {};
    \node[myv, fill=black] (v2) at (0.5,0) {};
    \node[myv, fill=black] (v3) at (1,0) {};
    \node[myv, fill=black] (v4) at (1.5,-0.25) {};
    \node[myv] (v5) at (0.75,-0.5) {};
    \draw (v1) -- (v2);
    \draw (v2) -- (v3);
    \draw (v3) -- (v4);
    \draw (v5) -- (v1);
    \draw (v5) -- (v2);
    \draw (v5) -- (v3);
    \draw (v5) -- (v4);
\end{tikzpicture} & 
\begin{tikzpicture}[myv/.style={circle, draw, inner sep=1.3pt}, scale=1]
    \node (o) at (0,0) {};
    \node[myv, fill=black] (v1) at (0,-0.25) {};
    \node[myv, fill=black] (v2) at (0.5,0) {};
    \node[myv] (v3) at (1,0) {};
    \node[myv, fill=black] (v4) at (1.5,-0.25) {};
    \node[myv] (v5) at (0.75,-0.5) {};
    \draw (v1) -- (v2);
    \draw (v2) -- (v3);
    \draw (v3) -- (v4);
    \draw (v5) -- (v1);
    \draw (v5) -- (v2);
    \draw (v5) -- (v3);
    \draw (v5) -- (v4);
\end{tikzpicture} & Gem & $3$ & $4$ \\  \cline{2-6}
 & \begin{tikzpicture}[myv/.style={circle, draw, inner sep=1.3pt}, scale=1]
    \node (o) at (0,0) {};
    \node[myv, fill=black] (v1) at (0,-0.25) {};
    \node[myv, fill=black] (v2) at (0.5,0) {};
    \node[myv, fill=black] (v3) at (1,0) {};
    \node[myv] (v4) at (1.5,-0.25) {};
    \node[myv, fill=black] (v5) at (0.75,-0.5) {};
    \draw (v1) -- (v2);
    \draw (v2) -- (v3);
    \draw (v3) -- (v4);
\end{tikzpicture} & 
\begin{tikzpicture}[myv/.style={circle, draw, inner sep=1.3pt}, scale=1]
    \node (o) at (0,0) {};
    \node[myv, fill=black] (v1) at (0,-0.25) {};
    \node[myv, fill=black] (v2) at (0.5,0) {};
    \node[myv, fill=black] (v3) at (1,0) {};
    \node[myv, fill=black] (v4) at (1.5,-0.25) {};
    \node[myv, fill=black] (v5) at (0.75,-0.5) {};
    \draw (v1) -- (v2);
    \draw (v2) -- (v3);
    \draw (v3) -- (v4);
\end{tikzpicture} & $\overline{\text{Gem}}$ & $5$ & $4$ \\ \hline \hline
\multirow{2}{*}{\X=\DTD} & \begin{tikzpicture}[myv/.style={circle, draw, inner sep=1.3pt}, scale=1]
    \node (o) at (0,0) {};
    \node[myv, fill=black] (v1) at (0,0.25) {};
    \node[myv, fill=black] (v2) at (0.5,0) {};
    \node[myv, fill=black] (v3) at (1,0) {};
    \node[myv, fill=black] (v4) at (1.5,0.25) {};
    \node[myv] (v5) at (0.75,-0.4) {};
    \draw (v1) -- (v2);
    \draw (v2) -- (v3);
    \draw (v3) -- (v4);
    \draw (v5) -- (v2);
    \draw (v5) -- (v3);
\end{tikzpicture} & 
\begin{tikzpicture}[myv/.style={circle, draw, inner sep=1.3pt}, scale=1]
    \node (o) at (0,0) {};
    \node[myv] (v1) at (0,0.25) {};
    \node[myv, fill=black] (v2) at (0.5,0) {};
    \node[myv, fill=black] (v3) at (1,0) {};
    \node[myv] (v4) at (1.5,0.25) {};
    \node[myv, fill=black] (v5) at (0.75,-0.4) {};
    \draw (v1) -- (v2);
    \draw (v2) -- (v3);
    \draw (v3) -- (v4);
    \draw (v5) -- (v2);
    \draw (v5) -- (v3);
\end{tikzpicture} & Bull & $3$ & $4$ \\ \cline{2-6}
 & 
 \begin{tikzpicture}[myv/.style={circle, draw, inner sep=1.3pt}, scale=1]
    \node (o) at (0,0) {};
    \node[myv] (v1) at (0,-0.25) {};
    \node[myv, fill=black] (v2) at (0.5,0) {};
    \node[myv, fill=black] (v3) at (1,0) {};
    \node[myv, fill=black] (v4) at (1.5,0) {};
    \node[myv] (v5) at (2,-0.25) {};
    \node[myv] (v6) at (1,-0.4) {};
    \draw (v1) -- (v2);
    \draw (v2) -- (v3);
    \draw (v3) -- (v4);
    \draw (v4) -- (v5);
    \draw (v6) -- (v3);
\end{tikzpicture} &
\begin{tikzpicture}[myv/.style={circle, draw, inner sep=1.3pt}, scale=1]
    \node (o) at (0,0) {};
    \node[myv, fill=black] (v1) at (0,-0.25) {};
    \node[myv, fill=black] (v2) at (0.5,0) {};
    \node[myv, fill=black] (v3) at (1,0) {};
    \node[myv, fill=black] (v4) at (1.5,0) {};
    \node[myv, fill=black] (v5) at (2,-0.25) {};
    \node[myv] (v6) at (1,-0.4) {};
    \draw (v1) -- (v2);
    \draw (v2) -- (v3);
    \draw (v3) -- (v4);
    \draw (v4) -- (v5);
    \draw (v6) -- (v3);
\end{tikzpicture} & Bow & $5$ & $3$ \\ \hline \hline
\multirow{2}{*}{\X=\LTD} & \begin{tikzpicture}[myv/.style={circle, draw, inner sep=1.3pt}, scale=1]
    \node (o) at (0,0) {};
    \node[myv, fill=black] (v1) at (0,0) {};
    \node[myv, fill=black] (v2) at (0,0.5) {};
    \node[myv, fill=black] (v3) at (0.5,0) {};
    \node[myv, fill=black] (v4) at (0.5,0.5) {};
    \draw (v1) -- (v2);
    \draw (v3) -- (v4);
\end{tikzpicture} &
\begin{tikzpicture}[myv/.style={circle, draw, inner sep=1.3pt}, scale=1]
    \node (o) at (0,0) {};
    \node[myv, fill=black] (v1) at (0,0) {};
    \node[myv] (v2) at (0,0.5) {};
    \node[myv, fill=black] (v3) at (0.5,0) {};
    \node[myv, fill=black] (v4) at (0.5,0.5) {};
    \draw (v1) -- (v2);
    \draw (v3) -- (v4);
\end{tikzpicture} & $2P_2$ & $3$ & $4$ \\ \cline{2-6}
 & \begin{tikzpicture}[myv/.style={circle, draw, inner sep=1.3pt}, scale=1]
    \node (o) at (0,0) {};
    \node[myv] (v1) at (0,0) {};
    \node[myv, fill=black] (v2) at (0.5,0) {};
    \node[myv, fill=black] (v3) at (1,0) {};
    \node[myv] (v4) at (1.5,0) {};
    \draw (v1) -- (v2);
    \draw (v2) -- (v3);
    \draw (v3) -- (v4);
\end{tikzpicture} &
\begin{tikzpicture}[myv/.style={circle, draw, inner sep=1.3pt}, scale=1]
    \node (o) at (0,0) {};
    \node[myv, fill=black] (v1) at (0,0) {};
    \node[myv, fill=black] (v2) at (0.5,0) {};
    \node[myv] (v3) at (1,0) {};
    \node[myv, fill=black] (v4) at (1.5,0) {};
    \draw (v1) -- (v2);
    \draw (v2) -- (v3);
    \draw (v3) -- (v4);
\end{tikzpicture} & $P_4$ & $3$ & $2$ \\ \hline
\end{tabular}
\caption{Comparison of \OSD-numbers and other \X-numbers of some graphs, where $\X \in \{\ID, \DTD, \LTD\}$. The black vertices constitute the respective codes.} \label{tab_comparisons}
\end{table}

\section{Computational complexity of the \OSD-problem} 
\label{sec3}

It has been established in the literature that all the previously studied X-problems are NP-hard~\cite{CHL2003,CS1987,SS2010}.
We next address the computational complexity of the decision version of the \OSD-problem\ defined as follows.

\defproblem{\sc \OD-Code}{An \OD-admissible graph $G$ and a positive integer $k$.}{Does there exist an \OD-code $S$ of $G$ such that $|S| \le k$?}

In this section, we shall prove that {\sc \OD-Code} is \NP-complete. Apart from that, we also consider in this section the computational complexity of another related problem as we now describe. By Theorem \ref{thm_osd_bounds-X}(b), given an \oldadmis ~graph $G$, we have either 
\[
\otd(G) = \od(G)~~ \text{ or }~~ \otd(G) = \od(G)+1.
\]
However, despite this small difference (of at most~$1$) between the two code-numbers, we show that it is \NP-hard to decide in general if they are equal or different on any given \OTD-admissible graph. This is equivalent to simply proving that it is \NP-hard to decide if $\otd(G) = \od(G) + 1$. We now define the decision version of this problem as the following.

\defproblem{\OTD\ = \OD+1}{An \OTD-admissible graph $G$ and an integer $k$.}
{Is $\otd(G) = k$ and $\od(G) = k-1?$ That is, are the following assertions true?
\begin{enumerate}[leftmargin=20pt, itemsep=0pt]
\item[(a)] There exists an \OTD-code $S$ of $G$ such that  $|S| = k$.
\item[(b)] There exists an \OD-code $S'$ of $G$ such that $|S'| = k-1$.
\item[(c)] For any vertex subset $S''$ of $G$, if $|S''| = |S|-1$, then $S''$ is not an \OTD-code of $G$; and if $|S''| = |S'|-1$, then $S''$ is not an \OD-code of $G$.
\end{enumerate}
}

As can be noticed in \OTD\ = \OD+1, given a vertex subset $S$ (with $|S| = k$) of an input graph $G$ on $n$ vertices, to check condition (c), one has to consider all subsets of $V(G)$ of order $k-1$ (which admits a running time of $\calO(n^{k-1})$) and check if $S$ is an \OTD-code (which admits a running time of $\calO(n^2)$). In other words, verifying a certificate for \OTD\ = \OD+1 can potentially take up to a running time of order $\calO(n^{k-1}) \cdot \calO(n^2) \subseteq \calO(n^{k+1})$. Since we have $k = |S| \in \calO(n)$, it implies that verifying a certificate for this problem may not necessarily be a polynomial-time algorithm. This implies that, \OTD\ = \OD+1 does not necessarily belong to the class \NP. Nevertheless, we show that the problem is still \NP-hard.

To prove the hardness of the above two decision problems, we use \textsc{Linear SAT}, or \textsc{LSAT} --- a variation of \textsc{3-SAT} --- which is defined as follows.

\defproblem{\textsc{Linear SAT (LSAT)}}{A set $X$ of Boolean variables and a set $\calC$ of clauses defined in terms of the literals from $X$ such that each clause contains at most $3$ literals; each literal from $X$ can appear in at most two clauses; and any two distinct clauses can have at most one literal in common.}{Does there exist a satisfying assignment on $X$?}

\textsc{LSAT} was proven to be \NP-complete by Arkin et al. in~\cite{ABCCKMS_2018}. From here on, we shall denote by $\psi = (X, \calC)$ a generic instance defined on the set of Boolean variables $X$ and the set of clauses $\calC$. In addition, let $n = |X|$ and $m = |\calC|$ and let us denote by $\phi = \phi(\calC)$ the Boolean formula defined by $\phi = \mathbf{c_1} \wedge \mathbf{c_2} \wedge \ldots \wedge \mathbf{c_m}$, where $\mathbf{c_i} \in \calC$ for all $1 \le i \le m$. We may also assume from here on that no clause contains both literals of a variable, as such a clause is always satisfiable.

We now introduce the following variation of LSAT where we require that if a literal from $X$ appears in a clause, it must appear in \emph{exactly} two clauses.

\defproblem{\textsc{Saturated Linear SAT (SL-SAT)}}{An instance $\psi = (X, \calC)$ of \text{LSAT} such that each literal from $X$ appears in either no clause or exactly two clauses.}{Does there exist a satisfying assignment on $X$?}

Let the \emph{auxiliary graph} of \textsc{SL-SAT} be the bipartite graph where one part of the vertex set is $\calC$ and the other the set of all literals from $X$. Moreover, let the edge set of the auxiliary graph consist of edges $x \mathbf{c}$, where $x$ is a literal from $X$ belonging to the clause $\mathbf{c} \in \calC$.

\begin{remark} \label{rem:SL-SAT}
The auxiliary graph of \textsc{SL-SAT} is bipartite, of maximum degree~$3$ and is of girth at least~$6$.
\end{remark}

\begin{proof}
Let $H$ be the auxiliary graph of \textsc{SL-SAT}. Then $H$ is bipartite by construction.

Each vertex of $H$ corresponding to the literals from $X$ has degree~$2$ (since a literal belongs to at most two clauses); and each vertex of $H$ corresponding to the clauses in $\calC$ has
degree at most~$3$ (since each clause contains at most~$3$ literals). This proves the result on the maximum degree of $H$.

Finally, $H$ cannot contain an induced four-cycle as it would imply that there exist two clauses having 
two literals in common which contradicts the input constraint of \textsc{SL-SAT}. Moreover, since $H$ is bipartite and hence, contains only even cycles, all induced cycles in $H$ must be of length at least~$6$. This proves the result on the girth of $H$.
\end{proof}

We shall next prove that \textsc{SL-SAT} is \NP-complete using the following reduction from an instance $\psi = (X,\calC)$ of \textsc{LSAT} to an instance $\psi' = (X', \calC')$ of \textsc{SL-SAT}.

\begin{reduction} \label{red:SL-SAT} \rm
Let $\phi$ be the Boolean formula of \textsc{LSAT} with respect to $X$ and $\calC$. Then the reduction starts with $X' = X$, $\calC' = \calC$ and $\phi' = \phi$. It then runs over all literals from $X$. At each step, if a literal, say $x$, appears in only one clause of $\phi$, the reduction creates an auxiliary variable $y$ and updates the set of variables $X' \leftarrow X' \cup \{y\}$, the set of clauses $\calC' \leftarrow \calC' \cup \{(x \vee y), (y)\}$ and also appends the formula $\phi \leftarrow (x \vee y) \wedge (y) \wedge \phi$. At the end of the reduction, the resulting $\psi' = (X', \calC')$ is an instance of \textsc{SL-SAT}. 
\end{reduction}

It can be verified that, at the end of Reduction~\ref{red:SL-SAT}, we have $|X'| \le 3n$ and $|\calC| \le m + 4n$. Moreover, the reduction is clearly polynomial-time in $n+m$.

\begin{lemma}
{\sc SL-SAT} is \NP-complete.
\end{lemma}

\begin{proof}
Clearly, \textsc{SL-SAT} belongs to the class \NP. Let $X' = X \cup Y$ and $\phi' = \phi \wedge \mu$, where $Y$ is the set of auxiliary
variables introduced in Reduction~\ref{red:SL-SAT} and $\mu$ is the Boolean formula obtained by appending the  clauses in Reduction~\ref{red:SL-SAT}. Note that $\phi$ has literals only from $X$. 

To prove the necessary part of the \NP-hardness, let $\psi = (X, \calC)$ be a \yes-instance of \textsc{LSAT} and so, let $A : X \to \{0,1\}$ be a satisfying assignment on $X$. Then, we create an assignment $A' : X' \to \{0,1\}$ by requiring $A'(x) = A(x)$ if $x \in X$ and $A'(x) = 1$ if $x \in Y$. Then it is clear that this is a satisfying assignment on $X'$ since $\phi' = \phi \wedge \mu$ and each clause in $\mu$ has an auxiliay variable which is assigned to~$1$ by $A'$. Hence, $\psi' = (X', \calC')$ is a \yes-instance of \textsc{SL-SAT}.

To prove the sufficiency part of the \NP-hardness, if $\psi'$ is a \yes-instance of \textsc{SL-SAT} and $A' : X' \to \{0,1\}$ is a satisfying assignment on $X'$, we create an assignment of $A : X \to \{0,1\}$ by requiring $A(x) = A'(x)$. Again, $\phi' = \phi \wedge \mu$ implies that $A$ is a satisfying assignment on $X$. Hence, $\psi$ is a \yes-instance of \textsc{LSAT}.
\end{proof}

We now come to the main results of this section and their proofs. To that end, we describe the following reduction from \textsc{SL-SAT} to both {\sc \OD-Code} and \OTD\ = \OD+1.

\begin{reduction} \label{red_NP_OD} \rm
The reduction takes as input an instance $\psi = (X,\calC)$ of \textsc{SL-SAT}. Also, let $\phi$ be the Boolean formula of \textsc{SL-SAT}. The reduction constructs a graph $G^\psi$ as follows (also refer to Figure~\ref{fig_NP_OD}):
\begin{itemize}[leftmargin=12pt, itemsep=0pt]
\item For every variable $x \in X$, do the following (refer to Figure~\ref{fig_NP Variable_OD}):
\begin{itemize}[leftmargin=8pt, itemsep=0pt]
\item Add a vertex named $w^x_1$ if the literal $x$ appears in the formula $\phi$; and add a vertex named  $w^x_2$ if the literal $\neg x$ appears in the formula $\phi$. For the rest of the reduction, we assume that both $w^x_1$ and $w^x_2$ exist; and in case one of them does not, the corresponding action related to the vertex described in this reduction does not take place.
\item Add 3 vertices named $v^x_1$, $v^x_2$ and $v^x_3$ and add edges $v^x_1 w^x_1$, $v^x_1 w^x_2$ and the edges $v^x_1 v^x_2$ and $v^x_2 v^x_3$ (the last two edges make the vertices $v^x_1$, $v^x_2$ and $v^x_3$ induce a $P_3$).
\end{itemize}
Let the graph induced by the above 5 (or 4) vertices be denoted as $G^x$ and be called the \emph{variable gadget} corresponding to the variable $x \in X$.
\item For every clause $\mathbf{c} \in \calC$, do the following (refer to Figure~\ref{fig_NP Clause_OD}):
\begin{itemize}[leftmargin=8pt, itemsep=0pt]
\item Add 3 vertices named $u^x_1$, $u^x_2$ and $u^x_3$ and add edges $u^x_1u^x_2$ and $u^x_2u^x_3$ thus making $u^x_1$, $u^x_2$ and $u^x_3$ induce a $P_3$.
\end{itemize}
Let the $P_3$ induced by the above 3 vertices be denoted by $P^\mathbf{c}$ and be called the \emph{clause gadget} corresponding to the clause $\mathbf{c} \in \calC$.
\item For all variables $x \in X$ and all clauses $\mathbf{c} \in \calC$, if the literal $x$ is in a clause $\mathbf{c}$, then add the edge $u^\mathbf{c}_1w^x_1$; and if the literal $\neg x$ is in a clause $\mathbf{c}$, then add the edge $u^\mathbf{c}_1w^x_2$ (refer to Figure~\ref{fig_NP G_I_OD}).
\end{itemize}
\end{reduction}

As a matter of notation, for any variable $x \in X$, let us denote by $[x]$ the set $\{x,\neg x\}$. Moreover, in all the proofs that follow in this section, we also assume that both the vertices $w^x_1$ and $w^x_2$, for $x \in X$, exist in the graph $G^\psi$. In case one of them does not and the analysis changes as a result of that, we shall point that out accordingly.
\begin{figure}[!t]
\centering
\begin{subfigure}[t]{0.3\textwidth}
\centering
\begin{tikzpicture}[blacknode/.style={circle, draw=black!, fill=black!, thick, scale=0.7},
whitenode/.style={circle, draw=black!, fill=white!, thick, scale=0.7},
scale=0.4]


\node[whitenode] (1) at (0,0) {};\node at (-1,0) {$w^x_1$};
\node[whitenode] (2) at (4,0) {};\node at (5,0) {$w^x_2$};
\node[whitenode] (3) at (2,2) {};\node at (3,2) {$v^x_1$};
\node[whitenode] (4) at (2,4) {};\node at (3,4) {$v^x_2$};
\node[whitenode] (5) at (2,6) {};\node at (3,6) {$v^x_3$};

\draw[-, thick, black!] (1) -- (3);
\draw[-, thick, black!] (2) -- (3);
\draw[-, thick, black!] (3) -- (4);
\draw[-, thick, black!] (4) -- (5);


\end{tikzpicture}
\caption{$G^x$: Variable gadget corresponding to $x \in X$.} \label{fig_NP Variable_OD}
\end{subfigure}
\hspace{2mm}
\begin{subfigure}[t]{0.3\textwidth}
\centering
\begin{tikzpicture}[blacknode/.style={circle, draw=black!, fill=black!, thick, scale=0.7},
whitenode/.style={circle, draw=black!, fill=white!, thick, scale=0.7},
scale=0.4]


\node[whitenode] (1) at (0,4) {};\node at (1,4) {$u^\mathbf{c}_1$};
\node[whitenode] (2) at (0,2) {};\node at (1,2) {$u^\mathbf{c}_2$};
\node[whitenode] (3) at (0,0) {};\node at (1,0) {$u^\mathbf{c}_3$};

\draw[-, thick, black!] (1) -- (2);
\draw[-, thick, black!] (2) -- (3);


\end{tikzpicture}
\caption{$P^\mathbf{c}$: Clause gadget corresponding to $\mathbf{c} \in \calC$.} \label{fig_NP Clause_OD}
\end{subfigure}
\hspace{2mm}
\begin{subfigure}[t]{0.3\textwidth}
\centering
\begin{tikzpicture}[blacknode/.style={circle, draw=black!, fill=black!, thick, scale=0.7},
whitenode/.style={circle, draw=black!, fill=white!, thick, scale=0.7},
scale=0.4]


\node[blacknode] (1) at (0,0) {};\node at (-1,0) {$w^x_1$};
\node[whitenode] (2) at (4,0) {};\node at (5,0) {$w^x_2$};
\node[blacknode, fill=black!30] (3) at (2,2) {};\node at (3,2) {$v^x_1$};
\node[blacknode] (4) at (2,4) {};\node at (3,4) {$v^x_2$};
\node[whitenode] (5) at (2,6) {};\node at (3,6) {$v^x_3$};

\node[blacknode] (6) at (-2,-2) {};\node at (-1,-2) {$u^\mathbf{c}_1$};
\node[blacknode] (7) at (-2,-4) {};\node at (-1,-4) {$u^\mathbf{c}_2$};
\node[whitenode] (8) at (-2,-6) {};\node at (-1,-6) {$u^\mathbf{c}_3$};

\draw[-, thick, black!] (1) -- (3);
\draw[-, thick, black!] (2) -- (3);
\draw[-, thick, black!] (3) -- (4);
\draw[-, thick, black!] (4) -- (5);

\draw[-, thick, black!] (6) -- (7);
\draw[-, thick, black!] (7) -- (8);

\draw[-, thick, black!] (6) -- (1);
\draw[dashed, thick, black!] (6) -- (-3,-0.5);
\draw[dashed, thick, black!] (6) -- (-4,-0.5);

\draw[dashed, thick, black!] (1) -- (1.7,-1.5);

\draw[dashed, thick, black!] (2) -- (6,-1.5);


\end{tikzpicture}
\caption{$G^\psi$: Instance of \mincode{\OD} and \OD~$\ne$ \OTD.} \label{fig_NP G_I_OD}
\end{subfigure}
\caption[Reduction~\ref{red_NP_OD} from \textsc{SL-SAT} to \mincode{\OD} and \OD~$\ne$ \OTD]{Polynomial-time construction of the graph $G^\psi$ from an \textsc{SL-SAT} ~instance $\psi = (X,\calC)$ as in Reduction~\ref{red_NP_OD}. The black vertices in (c) represent those in a code described in Lemma~\ref{lem_NP_1}. The gray vertex ($v^x_1$) implies that, for some fixed variable $x=x_0 \in X$, the vertex is not included in an \OD-code but is included in the \OTD-code described in Lemma~\ref{lem_NP_1}.}
\label{fig_NP_OD}
\end{figure}

\begin{lemma} \label{lem_NP_lb}
For an instance $\psi = (X,\calC)$ of \textsc{SL-SAT} with $|X|=n$ and $|\calC|=m$, let $G^\psi$ be as in Reduction~\ref{red_NP_OD}. Moreover, let
$T = \{(v^x_1, v^x_2, v^x_3) : x \in X\} \cup \{(u^\mathbf{c}_1,u^\mathbf{c}_2, u^\mathbf{c}_3) : \mathbf{c} \in \calC\}$.
Then, for an open-separating set $S$ of $G^\psi$, the following assertions are true.
 \begin{enumerate}[leftmargin=18pt, itemsep=0pt]
\item[(1)] If $S$ is an \OTD-code of $G^\psi$, then $|T \cap S| \ge 2n+2m$.
\item[(2)] If $S$ is an \OD-code of $G^\psi$, then $|T \cap S| \geq 2n+2m-1$.
\end{enumerate}
Moreover, we have $\otd(G^\psi) \geq 3n+2m$ and $\od(G^\psi) \geq 3n+2m-1$.
\end{lemma}

\begin{proof}
(1) Let $S$ be an \OTD-code of $G^\psi$. Then, for each $(v_1,v_2,v_3) \in T$, the vertex $v_2$ must belong to $S$ for the latter to total-dominate $v_3$. Moreover, at least one vertex from the pair $(v_1,v_3)$ must be in $S$ for the latter to total-dominate the vertex $v_2$. In total therefore, counting over all the $n$ variable gadgets and $m$ clause gadgets, we have $|T \cap S| \geq 2n+2m$.

\medskip

(2) Let $S$ be an \OD-code of $G^\psi$. Now, if at least two vertices out of each triple $(t_1,t_2,t_3) \in T$ belong to $S$, then we have $|S| \geq 2n+2m > 2n+2m-1$ and thus, we are done. So, let us assume that there exists one triple $(t^\star_1,t^\star_2,t^\star_3) \in T$ which has at most one vertex in $S$. However, every triple $(t_1,t_2,t_3) \in T$ must have at least one of $t_2$ and $t_3$ in $S$ in order for the latter to dominate $t_3$. This implies that exactly one of $t^\star_2$ and $t^\star_3$ belongs to $S$. This further implies that either $N_{G^\psi}(t^\star_2) \cap S = \emptyset$ or $N_{G^\psi}(t^\star_3) \cap S = \emptyset$. This further implies that each triple of $T$ other than $(t^\star_1,t^\star_2,t^\star_3)$ has at least two vertices in $S$ (or else, two vertices of $G^\psi$ would have empty neighborhoods in $S$ and hence, the pair would not be open-separated by $S$, a contradiction). Thus, again counting over all the $n$ variable gadgets and $m$ clause gadgets, we have $|S| \geq 2n+2m-1$.

Finally, in each of the above two cases, since $S$ is an open-separating set of $G^\psi$, at least one of $w^x_1$ and $w^x_2$ must be in $S$ in order for the latter to open-separate the pair $v^x_1, v^x_3$. This implies the final statement of the result.
\end{proof}
\begin{lemma} \label{lem_NP_1}
For an instance $\psi = (X,\calC)$ of \textsc{SL-SAT} with $|X| = n$ and $|\calC| = m$, let $G^\psi$ be the graph as constructed in Reduction~\ref{red_NP_OD}. Then, the existence of a satisfying assignment on $X$ implies that
\begin{enumerate}[leftmargin=18pt, itemsep=0pt]
\item[(1)] $\od(G^\psi) = 3n+2m-1$; and
\item[(2)] $\otd(G^\psi) = 3n+2m$.
\end{enumerate}
\end{lemma}

\begin{proof}
Let $x_0 \in X$ be any fixed variable. We now look at the two cases.

\medskip

(1) Since $\od(G^\psi) \geq 3n+2m-1$ by Lemma~\ref{lem_NP_lb}, in order to prove that $\od(G^\psi) = 3n+2m-1$, it is enough to show the existence of an \OD-code of $G^\psi$ of order at most $3n+2m-1$.  Let us now construct a code $S$ of $G^\psi$ by including in it the following vertices.
\begin{enumerate}[leftmargin=18pt, itemsep=0pt]
\item[(a)] $v^x_1$ for all $x \in X \setminus \{x_0\}$ and $v^x_2$ for all $x \in X$.
\item[(b)] $u^\mathbf{c}_1,u^\mathbf{c}_2$ for all $\mathbf{c} \in \calC$.
\item[(c)] For some $x \in X$, if exactly one of $x$ and $\neg x$ appears in the formula $\phi$, say $x$, then we pick $w^x_1$ in $S$, or else, we pick $w^x_2$ in $S$.
\item[(d)] If for some variable $x \in X$, both the literals $x$ and $\neg x$ belong to some respective clauses, say $\bf c$ and $\mathbf{c}'$, then exactly one of the literals is assigned~$1$ in the satisfying assignment on $X$. If $x$ is assigned~$1$, we pick $w^x_1$ in $S$, or else, we pick $w^x_2$ in $S$.
\end{enumerate}
This implies that every variable gadget has either $w^x_1$ or $w^x_2$ (but not both) in $S$ and hence, $|S| = 3n+2m-1$. We now show that $S$ is an \OD-code of $G^\psi$. To start with, we observe that $S$ is a dominating set of $G^\psi$. Thus, it is left to show that $S$ is also an open-separating set of $G$. We also observe that for exactly the one vertex $v^{x_0}_2$ of $G^\psi$, we have $N(v^{x_0}_2) \cap S = \emptyset$. Hence, to prove that $S$, indeed, is an open-separating set of $G^\psi$, by Remark~\ref{lem_char_osd}, it is enough to check that all pairs of vertices of $G^\psi$ of distance at most two between them are open-separated by $S$. We show next that this is true. Let $x \in X$ and $\bf c \in \calC$ be any general Boolean variable and clause, respectively.
\begin{itemize}[leftmargin=12pt, itemsep=0pt]
\item Since $N(v^{x_0}_2) \cap S = \emptyset$ uniquely, therefore, $S$ open-separates $v^{x_0}_2$ from every other vertex of $G^\psi$.

\item The vertex $v^x_2$ open-separates $v^x_3$ from every vertex of $G^\psi$ except $v^x_1$. However, the vertices $v^x_1, v^x_3$ are open-separated by whichever of $w^x_1$ and $w^x_2$ is in $S$.

\item Similarly, the vertex $u^\mathbf{c}_2$ open-separates $u^\mathbf{c}_3$ from every vertex of $G^\psi$ except $u^\mathbf{c}_1$. However, the clause $\mathbf{c}$ has a literal, say $x$ or $\neg x'$, which is assigned~$1$ in the satisfying assignment on $X$. Then, the vertices $u^\mathbf{c}_1, u^\mathbf{c}_3$ are open-separated by $w^x_1$ if $x$ is assigned~$1$ or by $w^x_2$ if $\neg x'$ is assigned~$1$.

\item Every vertex of $G^\psi$ in $S$ open-separates itself from all its neighbors. Therefore, we now check that each vertex in $S$ is also open-separated by $S$ from all other vertices at distance~$2$ from the former.

\item The vertex $v^x_2$ open-separates $v^x_1$ from every vertex of $G^\psi$ except $v^x_3$. However, the pair $v^x_1, v^x_3$ were previously shown to be open-separated by $S$.

\item Similarly, the vertex $u^\mathbf{c}_2$ open-separates $u^\mathbf{c}_1$ from every vertex of $G^\psi$ except $u^\mathbf{c}_3$. Once again, the pair $u^\mathbf{c}_1, u^\mathbf{c}_3$ was already shown to be open-separated by $S$.

\item We now look at the open-separation of the pair $v^x_2, w^x_i$ by $S$, where $i \in [2]$. Let us assume, without loss of generality, that the vertex $w^x_1$ exists in the graph $G^\psi$. Then, $w^x_1$ has a neighbor $u^\mathbf{c}_1$ for some clause $\mathbf{c} \in \calC$. Then, we have $w^x_1 u^\mathbf{c}_1 \in E(G^\psi)$. Since $u^\mathbf{c}_1 \in S$, the pair $v^x_2, w^x_1$ is open-separated by $u^\mathbf{c}_1 \in S$. The case for the open-separation of the pair $v^x_2, w^x_2$ (if $w^x_2$ exists in $G^\psi$) by $S$ follows by exactly the same reasoning.

\item We now look at the open-separation by $S$ of the vertex $u^\mathbf{c}_2$ from other vertices (of distance~$2$) from it. Let $x \in X$ be such that $x$ has a literal in $\mathbf{c}$. Without loss of generality, let $x$ itself be a literal in $\mathbf{c}$. Thus, we look at the pair $w^x_1, u^\mathbf{c}_2$ to be open-separated by $S$. Since the literal $x$ appears in two clauses, say $\mathbf{c}$ and $\mathbf{c}'$, we have the edges $w^x_1 u^\mathbf{c}_1 \in E(G^\psi)$ and $w^x_1 u^{\mathbf{c}'}_1 \in E(G^\psi)$. This implies that the vertices $w^x_1, u^\mathbf{c}_2$ are open-separated by $u^{\mathbf{c}'}_1$. Similarly, the vertices $w^x_1, u^{\mathbf{c}'}_2$ are open-separated by $u^\mathbf{c}_1$ (this is where we use the property of saturating an \textsc{LSAT} formula so that the literal $x$ appears in two clauses. Otherwise, if $x_0$ belongs to the clause $\mathbf{c_0}$, the pair $u^{\mathbf{c_0}}_2, w^{x_0}_1$ would not be open-separated by $S$ since $v^{x_0}_1, u^{\mathbf{c_0}}_3 \notin S$).
\end{itemize}

Finally, the only vertices in different variable gadgets (respectively, clause gadgets) with distance 2 between them are of the form $w^x_i \in G^x$ and $w^{x'}_j \in G^{x'}$ (respectively, $u^\mathbf{c}_1 \in P^\mathbf{c}$ and $u^{\mathbf{c}'}_1 \in P^{\mathbf{c}'}$), where $x,x' \in X$ (respectively, $\mathbf{c},\mathbf{c}' \in \calC$) are distinct. However, the set $S$ open-separates any such pair $w^x_i, w^{x'}_j$ by either $v^x_1$ (if $x_1 \ne x_0$) or by $v^{x'}_1$ (if $x'_1 \ne x_0$); and the pair $u^\mathbf{c}_1, u^{\mathbf{c}'}_1$ by $u^\mathbf{c}_2$. This proves that $S$ is an open-separating set of $G^\psi$.\\

(2) Again, by Lemma~\ref{lem_NP_lb}, we have $\otd(G^\psi) \geq 3n+2m$. Therefore, to show that $\otd(G^\psi) = 3n+2m$, it is enough to show that there exists an \OTD-code of $G^\psi$ of order at most $3n+2m$. To that end, we simply replace the \OD-code $S$ constructed above in part (1) by $S' = S \cup \{v^{x_0}_1\}$ which then becomes a total-dominating set of $G^\psi$. Moreover, since adding vertices to an open-separating set keeps it open-separating, $S'$ is also an open-separating set of $G^\psi$. Moreover, $|S'| = |S| + 1 = 3n+2m$. This proves the lemma.
\end{proof}

\begin{lemma} \label{aplem_NP_w_sat_assign_OD}
Let $\psi = (X,\calC)$ be an instance of \textsc{SL-SAT} and $G^\psi$ be the graph as in Reduction~\ref{red_NP_OD}. If there exists an open-separating set $S$ of $G^\psi$ such that $|\{w^x_1,w^x_2\} \cap S| = 1$ for all $x \in X$, then $X$ has a satisfying assignment.
\end{lemma}

\begin{proof}
Let $S$ be an open-separating set of $G^\psi$ such that $|\{w^x_1,w^x_2\} \cap S| = 1$ for all $x \in X$. We now provide a binary assignment on $X$ the following way: for any $x \in X$, if $w^x_1 \in S$ and $w^x_2 \notin S$, then put $(x,\neg x)= (1,0)$; and if $w^x_1 \notin S$ and $w^x_2 \in S$, then put $(x,\neg x) = (0,1)$. Clearly, this assignment on $X$ is a valid one since, for each $x \in X$, exactly one of $x$ and $\neg x$ is assigned $1$ and the other $0$. To now prove that this assignment on $X$ is also a satisfying one, we simply note that, for each clause $\mathbf{c} \in \calC$, in order for the set $S$ to open-separate $u^\mathbf{c}_1$ and $u^\mathbf{c}_3$, either there exists a variable $x \in X$ such that $u^\mathbf{c}_1w^x_1 \in E(G^\psi)$ with $w^x_1 \in S$; or there exists a variable $x' \in X$ such that $u^\mathbf{c}_1w^{x'}_2 \in E(G^\psi)$ and $w^{x'}_2 \in S$. Therefore, by construction of the graph $G^\psi$, either $x$ is a variable in the clause $\mathbf{c}$ with $(x,\neg x) = (1,0)$ or $x'$ is a variable whose literal $\neg x'$ is in the clause $\mathbf{c}$ with $(x', \neg x') = (0,1)$. Hence, the binary assignment formulated on $X$ is a satisfying one.
\end{proof}

\begin{lemma} \label{lem_NP_2}
For an instance $\psi$ of \textsc{SL-SAT} with $|X| = n$ and $|\calC| = m$, let $G^\psi$ be the graph as constructed in Reduction~\ref{red_NP_OD}. Then, the following assertions are true.
\begin{enumerate}[leftmargin=18pt, itemsep=0pt]
\item[(1)]\label{aplem_NP_2_OD} If $\od(G^\psi) = 3n+2m-1$, then $X$ has a satisfying assignment.
\item[(2)]\label{aplem_NP_2_OTD} If $\otd(G^\psi) = 3n+2m$, then $X$ has a satisfying assignment.
\end{enumerate}
\end{lemma}
\begin{proof}
Let $S$ be either an \OD-code or an \OTD-code of $G^\psi$. Since $S$ open-separates the vertices $v^w_1, v^x_3$, it implies that at least one of $w^x_1$ and $w^x_2$ must be in $S$. Then, by Lemma~\ref{aplem_NP_w_sat_assign_OD}, to show that $X$ has a valid satisfying assignment, it is enough to show that $|\{w^x_1, w^x_2\} \cap S| = 1$.

Let us first assume that $\od(G^\psi) = 3n+2m-1$ and let $S'$ be an \OD-code of $G^\psi$ such that $|S'| = 3n+2m-1$. If on the contrary, there exists some variable $x' \in X$ for which $|\{w^x_1, w^x_2\} \cap S'| = 2$, then using Lemma~\ref{lem_NP_lb}(1), we have
\begin{align*}
3n+2m-1 = |S'| &= |T \cap S'| + |\{w^{x'}_1, w^{x'}_2\} \cap S'| + \sum_{x \in X \setminus \{x'\}} |\{w^x_1, w^x_2\} \cap S'| \ge 3n+2m.
\end{align*}
which is a contradiction. Therefore, $|\{w^{x'}_1, w^{x'}_2\} \cap S'| = 1$ and this proves (1).

Let us now assume that $\otd(G^\psi) = 3n+2m$ and let $S$ be an \OTD-code of $G^\psi$ such that $|S| = 3n+2m$. Again, if on the contrary, there exists some variable $x' \in X$ for which $|\{w^{x'}_1, w^{x'}_2\} \cap S| = 2$, then using Lemma~\ref{lem_NP_lb}(2), we have
\begin{align*}
3n+2m = |S| &= |T \cap S| + |\{w^{x'}_1, w^{x'}_2\} \cap S| + \sum_{x \in X \setminus \{x'\}} |\{w^x_1, w^x_2\} \cap S| \ge 3n+2m+1.
\end{align*}
which is again a contradiction. Therefore, $|\{w^{x'}_1, w^{x'}_2\} \cap S| = 1$ and this proves (2).
\end{proof}

\begin{lemma} \label{lem:NP_OD_bibpartite}
For an instance $\psi$ of \textsc{SL-SAT}, let $G^\psi$ be the graph as constructed in Reduction~\ref{red_NP_OD}. Then, $G^\psi$ is bipartite and of maximum degree~$4$ and girth at least~$6$.
\end{lemma}

\begin{proof}
Let $V_1 = \{u^\mathbf{c}_1, u^\mathbf{c}_3, v^x_1, v^x_3 : x \in X, \mathbf{c} \in \calC \}$ and $V_2 = \{u^\mathbf{c}_2, w^x_1, w^x_2, v^x_2 : x \in X, \mathbf{c} \in \calC\}$. Then it can be checked that the sets $V_1$ and $V_2$ are each independent sets. Moreover, $V(G^\psi) = V_1 \cup V_2$. This implies that $G^\psi$ is a bipartite graph.

All vertices in $\{u^\mathbf{c}_2, u^\mathbf{c}_3, v^x_2, v^x_3\}$ have degree at most~$3$ and the vertex $v^x_1$ is of degree~$3$. Since each literal from $X$ can appear in at most two clauses, it implies that the degrees of the vertices $w^x_1$ and $w^x_2$ are at mos~$3$. Finally, every clause $\mathbf{c} \in \calC$ has at most~$3$ literals from $X$. This implies that the vertex $u^\mathbf{c}_1$ is of degree at most~$4$. This proves that the graph $G^\psi$ is of maximum degree~$4$.

Notice that any induced cycle of $G^\psi$ must include vertices only in the set $S = \{u^{\mathbf{c}}_1 : \mathbf{c} \in \calC\} \cup \{v^x_1, w^x_1, w^x_2 : x \in X\}$. Also, let $S' = \{u^{\mathbf{c}}_1 : \mathbf{c} \in \calC\} \cup \{w^x_1, w^x_2 : x \in X\}$. Then $G^\psi[S']$ is isomorphic to the auxiliary graph of \textsc{SL-SAT}. Therefore, by Remark~\ref{rem:SL-SAT}, any cycle of $G^\psi$ induced only by the vertices of $S'$ must be of length at least~$6$. Moreover, any cycle of $G^\psi$ with a vertex in $\{v^x_1 : x \in X\}$ must also include the vertices $w^x_1$ and $w^x_2$ with $G^\psi[v^x_1, w^x_1, w^x_3]$ being isomorphic to $P_3$. Since the latter two vertices are not adjacent to a common $u^{\mathbf{c}}_1$ (or else, both literals of a variable would appear in the same clause), the cycle must be of length at least~$6$. Hence, this proves that all induced cycles of $G^\psi$ are of length at least~$6$. Therefore, the graph $G^\psi$ is of girth at least~$6$.
\end{proof}

This brings us to the proofs of our main theorems in this section.

\begin{theorem} \label{thm_osd_hard}
\mincode{\OD} is \NP-complete, even when the input graph is bipartite and of maximum degree~$4$ and girth at least~$6$.
\end{theorem}

\begin{proof}
\mincode{\OD} clearly belongs to the class \NP\ since it can be verified in polynomial-time if a given vertex subset of a graph $G$ is an \OD-code of $G$. To prove \NP-hardness, we show that an instance $\psi=(X,\calC)$ with $|X|=n$ and $|\calC|=m$ is a \yes-instance of \textsc{SL-SAT} if and only if $(G^\psi,3n+2m-1)$ is a \yes-instance of \mincode{\OD}, where the graph $G^\psi$ is as constructed in Reduction~\ref{red_NP_OD}. In other words, we show that $X$ has a satisfying assignment if and only if there exists an \OD-code $S$ of $G^\psi$ such that $|S| \leq 3n+2m-1$.

To prove the necessary part of the last statement, if $X$ has a satisfying assignment, then by Lemma~\ref{lem_NP_1}(1), we have $\od(G^\psi) = 3n+2m-1$. This implies that there exists an \OD-code $S$ of $G^\psi$ such that $|S| = 3n+2m-1$. On the other hand, to prove the sufficiency, if there exists an \OD-code $S$ of $G^\psi$ such that $|S| \le 3n+2m-1$, then combining with Lemma~\ref{lem_NP_lb}, we have $\od(G^\psi) = 3n+2m-1$. Therefore, by Lemma~\ref{lem_NP_2}(1), there exists a satisfying assignment on $X$.

Finally, by Lemma~\ref{lem:NP_OD_bibpartite}, the graph $G^\psi$ is bipartite and of maximum degree~$4$ and girth at least~$6$. This proves the result.
\end{proof}
\begin{theorem} \label{thm_osd-old_hard}
\textsc{\OTD\ = \OD+1} is \NP-hard, even when the input graph is bipartite and of maximum degree~$4$ and girth at least~$6$.
\end{theorem}

\begin{proof}
We show that \OTD\ = \OD+1 is \NP-hard by showing that an instance $\psi=(X,\calC)$ with $|X|=n$ and $|\calC|=m$ is a \yes-instance of \textsc{SL-SAT} if and only if $(G^\psi, 3n+2m)$ is a \yes-instance of \OTD\ = \OD+1, where the graph $G^\psi$ is as constructed in Reduction~\ref{red_NP_OD}. In other words, we show that $X$ has a satisfying assignment if and only if $\od(G^\psi) = 3n+2m-1$ and $\otd(G^\psi) = 3n+2m$.

The necessary condition is true since, if $X$ has a satisfying assignment, then by Lemma~\ref{lem_NP_1}, we have $\otd(G^\psi) = 3n+2m$ and $\od(G^\psi) = 3m+2m-1$. To prove the sufficiency condition, let us assume that $\otd(G^\psi) = 3n+2m$ and $\od(G^\psi) = 3n+2m-1$. Then, by Lemma~\ref{lem_NP_2}, $X$ has a satisfying assignment.

Finally, by Lemma~\ref{lem:NP_OD_bibpartite}, the graph $G^\psi$ is bipartite and of maximum degree~$4$ and girth at least~$6$. This proves the result.
\end{proof}

Reduction~\ref{red_NP_OD} also proves \mincode{\OTD} to be \NP-complete. Even though this result was also shown in~\cite{SS2010} by Seo and Slater, we believe that our reduction is perhaps simpler than the one used in~\cite{SS2010}. This actually answers a
question posed by Slater~\cite{FS_2013} regarding the existence of an easier proof for the \NP-completeness of \mincode{\OTD}. Moreover, our result proves \mincode{\OTD} to be \NP-complete even for bipartite graphs of maximum degree~$4$ and girth at least~$6$ (thus, this is a stronger result with an arguably simpler construction than the one in~\cite{PP_2017} where the authors show the problem to be \NP-complete on bipartite graphs). We now present here the following formal proof.

\begin{theorem}
\mincode{\OTD} is \NP-complete even when the input graph is bipartite and of maximum degree~$4$ and girth at least~$6$.
\end{theorem}

\begin{proof}
\mincode{\OTD} clearly belongs to the class \NP\ since it can be verified in polynomial-time if a given vertex subset of a graph $G$ is an \OTD-code of $G$. To prove \NP-hardness, we show that an instance $\psi=(X,\calC)$ with $|X|=n$ and $|\calC|=m$ is a \yes-instance of \textsc{SL-SAT} if and only if $(G^\psi,3n+2m)$ is a \yes-instance of \mincode{\OTD}, where the graph $G^\psi$ is as constructed in Reduction~\ref{red_NP_OD}. In other words, we show that $X$ has a satisfying assignment if and only if there exists an \OTD-code $S$ of $G^\psi$ such that $|S| \leq 3n+2m$.

To prove the necessary part of the last statement, if $X$ has a satisfying assignment, then by Lemma~\ref{lem_NP_1}(2), we have $\od(G^\psi) = 3n+2m$. This implies that there exists an \OTD-code $S$ of $G^\psi$ such that $|S| = 3n+2m$. On the other hand, to prove the sufficiency, if there exists an \OTD-code $S$ of $G^\psi$ such that $|S| \le 3n+2m$, then combining with Lemma~\ref{lem_NP_lb}, we have $\otd(G^\psi) = 3n+2m$. Therefore, by Lemma~\ref{lem_NP_2}(2), there exists a satisfying assignment on $X$.

Finally, by Lemma~\ref{lem:NP_OD_bibpartite}, the graph $G^\psi$ is bipartite and of maximum degree~$4$ and girth at least~$6$. This proves the result.
\end{proof}

\section{OD-numbers of some graph families}
\label{sec4}
In this section, we study the \OD-numbers of graphs belonging to some well-known graph families. Moreover, 
motivated by the hardness of deciding for which graphs the \OD- and the \OTD-numbers differ, we compare in the following the two numbers on some chosen graph families. 
This comparison also exhibits extremal cases for the upper bounds in Theorem~\ref{thm_osd_bounds-n}. 
To do so, we shall first determine the \OD-hypergraphs of the studied graphs, and then deduce their \OD-numbers.

\subsection{\OD-Clutters}
For some previously studied \X-problems, it was observed, for instance, in \cite{ABLW_2018,ABLW_2022} that $\hyp_\X(G) = (V,\F_\X)$ may contain redundant hyperedges. 
In fact, if there are two hyperedges $F, F' \in \F_\X$ with $F \subset F'$, then $F \cap C \neq \emptyset$ also implies $F' \cap C \neq \emptyset$ for every $C \subset V$.
Thus, $F'$ is \emph{redundant} as $(V,\F_\X -\{F'\})$ suffices to encode the domination and separation properties of the X-codes of $G$. This motivates to consider the \emph{\X-clutter} $\mathcal{C}_\X(G)$ of the graph $G$ obtained from $\hyp_\X(G)$ by removing all redundant hyperedges of the latter, as clearly $\tau(\hyp_{\X}(G)) = \tau(\clu_{\X}(G))$ holds. Moreover, since isolated vertices of a hypergraph are not relevant for determining its covering number, we also assume from here on that $\calC_X(G)$ is without isolated vertices. More precisely, we therefore define 
\[\clu_{\X}(G) = (V \setminus V^0_{\X}(G),\F_{\X}^\star(G))
\]
to be the $\X$-clutter of $G$, where $\F_{\X}^\star(G)$ is the set of all non-redundant hyperedges of $\hyp_\X(G)$ and $V^0_{\X}(G)$ denotes the subset of all vertices of $G$ not belonging to any hyperedge in $\F_{\X}^\star(G)$.

Moreover, a special interest lies in hyperedges of $\mathcal{C}_\X(G)$ consisting of a single vertex, called a \emph{forced vertex}, as each forced vertex needs to belong to every X-code of $G$.
We denote the set of singleton hyperedges containing forced vertices of $G$ by $\for_{\X}(G)$ and the set $\F_\X^\star(G) \setminus \for_\X(G)$ of hyperedges of cardinality at least~$2$ by $\fac_\X(G)$. In other words, $\F_\X^\star(G) = \for_\X(G) \cup \fac_\X(G)$. 
For \OD-codes, observe that $\for_{\OD}(G)$ combines vertices forced by domination and vertices forced by open-separation. 
Accordingly, we deduce from \cite{ABLW_2022} the following:
\begin{corollary} \label{lem_osd_forced}
For an \odadmis ~graph $G$, we have $\for_{\OD}(G) = V_D \cup V_O$, where
\[
\begin{array}{ccl}
V_D & = & \{x \in V : N(x) = \emptyset\}; \text{ and} \\
V_O & = & \{y \in V : \exists \mbox{ non-adjacent } u,v \mbox{ with } \{y\} = N(u) \Delta N(v)\}. \\
\end{array}
\]
\end{corollary}

Note that for an \odadmis ~graph $G$ having an isolated vertex $v$, $V_O$ contains the support vertex of each leaf, say $u$, of $G$ (as $N(u) \Delta N(v) = N(u) \Delta \emptyset = N(u)$ holds). 

\begin{example} \rm
For illustration, we construct $\hyp_{\OD}(P_4)$ and $\mathcal{C}_{\OD}(P_4)$. 
The \OD-hypergraph $\hyp_{\OD}(P_4)$ is composed of 
\[
\begin{array}{lclcl}
N[1]=\{1,2\}   & \ \ & N(1) \Delta N(2) = \{1,2,3\}   & \ \ & N(1) \Delta N(3) = \{4\}\\
N[2]=\{1,2,3\} & & N(2) \Delta N(3) = \{1,2,3,4\} & & N(1) \Delta N(4) = \{2,3\}\\
N[3]=\{2,3,4\} & & N(3) \Delta N(4) = \{2,3,4\}   & & N(2) \Delta N(4) = \{1\}\\
N[4]=\{3,4\}   & & & &\\
\end{array}
\]
Clearly, the \OD-clutter $\mathcal{C}_{\OD}(P_4)$ only contains the symmetric differences of open neighborhoods of non-adjacent vertices, namely, the sets $\{1\}, \{2,3\}, \{4\}$. Moreover, we have $\for_{\OD}(P_4) = \{1,4\}$, $\fac_{\OD}(P_4) = \{\{2,3\}\}$ and $V^0_{\OD}(P_4) = \emptyset$.
\end{example}

Note that for the previously studied X-problems with $\X \in \{\ID,\ITD,\LD,\LTD,\OTD\}$, it has been shown in \cite{ABLW_2018,ABLW_2022} that $\mathcal{C}_\X(G)$ does not contain symmetric differences of neighborhoods of non-adjacent vertices without common neighbor. This does not apply to \OD-clutters, as $N(1) \Delta N(4) = \{2,3\}$ from the above example demonstrates.

The proofs in the remainder of this section rely on arguments from the literature and on arguments issued from the representation of the \OD-problem as covering problem in the \OD-clutter. For the latter, we use the notation $\Delta(u,v)= N(u)\bigtriangleup N(v)$ to refer to symmetric differences. 

In the sequel, it will turn out that the \OD-clutters of several studied graph families are related to the hypergraph ${\cal{R}}_n^q=(V,{\cal E})$ called \emph{complete $q$-rose of order $n$}, where $V=\{1,\ldots,n\}$ and $\cal E$ contains all $q$-element subsets of $V$ for $2 \leq q < n$.  
In \cite{S_1989} it was proved that we have $\tau({\cal{R}}_n^q)=n-q+1$.
Note that, for $q=2$, ${\cal{R}}_n^q$ is in fact the clique $K_n$ and $\tau(K_n)=n-1$.

\subsection{Cliques and related families}
Cliques $K_n$ are clearly open twin-free so that for $n \geq 2$ both \OD- and \OTD-codes exist.
It is known from e.g. \cite{ABLW_2022} that the following holds:
$$
\old(K_n) = \begin{cases}
            2, &\text{if } n=2; \\
            n-1, &\text{if } n \geq 3.
            \end{cases}
$$
\begin{lemma} \label{lem_OD_cliques}
For a clique $K_n$ with $n \geq 2$, we have $\mathcal{C}_{\OSD}(G)= {\cal{R}}_n^2=K_n$ and $\osd(K_n) = n-1$.
\end{lemma}

\begin{proof}
Consider a clique $K_n = (V, V \times V)$ with $n \geq 2$. 
$\hyp_{\OD}(K_n)$ is clearly composed of 
\begin{itemize}[leftmargin=12pt, itemsep=0pt]
\item the closed neighborhoods $N[v]=V$ of all vertices $v \in V$ and 
\item the symmetric differences $\Delta(u,v)= \{u,v\}$ of open neighborhoods of distinct vertices $u,v \in V$.
\end{itemize}
This shows that all neighborhoods are redundant. 
Thus, $\clu_{\OD}(K_n) = {\cal{R}}_n^2 = K_n$ follows which implies $\osd(K_n) = \tau({\cal{R}}_n^2) = n-1$ by \cite{S_1989}.
\end{proof}

Hence, the \OD- and the \OTD-numbers of cliques $K_n$ differ only for $n=2$ and are equal for all $n \geq 3$. Moreover, the upper bound in Theorem \ref{thm_osd_bounds-n} is attained for all $n \geq 3$.

Consider now a graph $G = K_{n_1} \oplus \ldots \oplus K_{n_k}$ that is the disjoint union of $k \geq 2$ cliques with $1 < n_1 \leq \ldots \leq n_k$.
It is well-known that the \OTD-number of the disjoint union of two or more graphs is the sum of their \OTD-numbers. Hence, we have
$$\otd(K_{n_1} \oplus \ldots \oplus K_{n_k}) = \sum_{n_i = 2} 2 + \sum_{n_i\geq 3} (n_i-1).$$
To compare this with the corresponding \OD-numbers, we have the following.
\begin{lemma} \label{lem_OD_union-cliques}
  Let $G = K_{n_1} \oplus \ldots \oplus K_{n_k}$ be a disjoint union of $k \geq 2$ cliques with $1 < n_1 \leq \ldots \leq n_k$.
  \begin{itemize}[leftmargin=18pt, itemsep=0pt]
  \item[(a)] If $n_1 = 2$, then $\od(G) = -1 + \sum_{n_i = 2} 2 + \sum_{n_i\geq 3} (n_i-1)$,
  \item[(b)] If $n_1\geq 3$, then $\od(G) = \sum_{1 \leq i \leq k} (n_i-1)$.
\end{itemize}
\end{lemma}

\begin{proof}
Consider the disjoint union $G = K_{n_1} \oplus \ldots \oplus K_{n_k}$ of cliques with $1 < n_1 \leq \ldots \leq n_k$ and $k \geq 2$ and suppose that $n_i = 2$ for all $i \leq \ell \leq k$ (with possibly $0 = \ell$, that is all $n_i \geq 3$).
Let us further denote $V_i = V(K_{n_i})$. Then, $\hyp_{\OD}(G)$ is clearly composed of 
\begin{itemize}[leftmargin=12pt, itemsep=0pt]
\item the closed neighborhoods $N[v]=V_i$ of all vertices $v \in V_i$, $1 \leq i \leq k$ and 
\item the symmetric differences $\Delta(u,v)= \{u,v\}$ of distinct vertices $u,v \in V_i$ as well as $\Delta(u,v)= N(u) \cup N(v)$ for $v \in V_i$, $u \in V_j$ with $i \neq j$.
\end{itemize}
This shows that all neighborhoods are redundant (by $1 < n_1$), as well as all symmetric differences of open neighborhoods of vertices from different components $K_{n_i}$ and $K_{n_j}$ as soon as at least one of $n_i$ and $n_j$ is $\ge 3$. This implies $\clu_{\OD}(G) = K_{2\ell} \oplus K_{n_{\ell+1}} \oplus \ldots + K_{n_k}$ and, accordingly, 
$\od(G) = (2\ell-1) + \sum_{\ell < i \leq k} (n_i-1)$ follows by \cite{S_1989}.
\end{proof}

Hence, for graphs $G$ that are disjoint unions of cliques, the \OD- and the \OTD-numbers are equal if all components are cliques of order $\geq 3$, but differ otherwise.
In particular, if $G = k K_2$ is a \emph{matching} (i.e. if $n_i = 2$ for all $1 \leq i \leq k$) and $k \geq 2$,
then the arguments from the previous proof show that $\clu_{\OD}(k K_2)=K_{2k}$ and thus $\od(k K_2)=2k-1$ holds. 
Thus, the \OTD-number of $G = k K_2$ is strictly greater than its \OD-number,
and the upper bound of $\od(G) = |V(G)| - 1$ from Theorem \ref{thm_osd_bounds-n} is attained.

\begin{figure}[h]
\begin{center}
\includegraphics[scale=1.0]{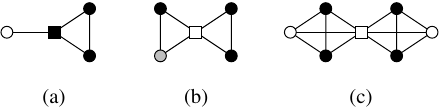}
\caption{Three $2$-clique-stars (the square indicates the universal vertex, black vertices belong to a minimum \OD-code, the grey vertex needs to be added to obtain an \OTD-code), the graph in (b) is a 2-fan.}
\label{fig_clique-stars}
\end{center}
\end{figure}

Next, we study graphs obtained from the disjoint union of cliques $K_{n_1} \oplus \ldots \oplus K_{n_k}$ for some $k \geq 2$ by adding a universal vertex and call them \emph{$k$-clique-stars} $S_k$, see Figure \ref{fig_clique-stars} for illustration. 
To have no open twins, at most one of the cliques can be of order one, that is, assuming $n_1 \leq \ldots \leq n_k$, we require that $n_2 \geq 2$ holds. 

\begin{theorem} \label{thm_OD_clique-stars}
  Consider a $k$-clique-star $S_k$ obtained from the disjoint union of cliques $K_{n_1} \oplus \ldots \oplus K_{n_k}$ for some $k \geq 2$ by adding a universal vertex and assume $n_1 \leq \ldots \leq n_k$.
  \begin{itemize}
  \item[(a)] If $n_1 = 1$ and $n_2 = \ldots = n_{\ell} = 2$ for some $\ell \leq k$, then we have
    $$\clu_{\OD}(S_k) = K_2 \oplus 2(\ell-1)K_1 \oplus \sum_{n_i\geq 3} K_{n_i} \mbox{ and } 
    \clu_{\OTD}(S_k) = (2\ell-1)K_1 \oplus \sum_{n_i\geq 3} K_{n_i}$$
    so that $\od(S_k) = \otd(S_k) = 2\ell-1 + \sum_{n_i\geq 3} (n_i-1)$ follows.
  \item[(b)] If $n_1 = \ldots = n_{\ell} = 2$ for some $1 \leq \ell \leq k$, then we have
        $$\clu_{\OD}(S_k) = K_{2\ell} \oplus \sum_{n_i\geq 3} K_{n_i} \mbox{ and } \clu_{\OTD}(S_k) = K_{2\ell+1} \oplus \sum_{n_i\geq 3} K_{n_i}$$
       so that $\od(S_k) =  2\ell-1 + \sum_{n_i\geq 3} (n_i-1)$ but $\otd(S_k) = 2\ell + \sum_{n_i\geq 3} (n_i-1)$ follows.
  \item[(c)]  If $n_1\geq 3$, then we have $\clu_{\OD}(S_k) = \clu_{\OTD}(S_k) = \sum_{1 \leq i \leq k} K_{n_i}$ and $\od(S_k) = \otd(S_k) = \sum_{1 \leq i \leq k} (n_i-1)$ follows.
\end{itemize}
\end{theorem}

\begin{proof}
  Consider such a $k$-clique-star $S_k$ for some $k \geq 2$ and let $V = \{u\} \cup \bigcup_{1 \leq i \leq k} V_i$ where $u$ is the universal vertex and $V_i$ denotes the vertex subset of $K_{n_i}$. Furthermore, let $W_2 = \bigcup_{n_i = 2} V_i$.
  In order to determine $\clu_{\X}(S_k)$ and $\x(S_k)$ for $\X \in \{\OD,\OTD\}$, we first construct the \X-hypergraphs.
  Recall that the following hyperedges are involved: the closed or open neighborhoods
\begin{itemize}[leftmargin=12pt, itemsep=0pt]
\item $N[v]=V_i \cup \{u\}$ and $N(v) = (V_i\setminus \{v\}) \cup \{u\}$ of all vertices $v \in V_i$, $1 \leq i \leq k$ and 
  \item $N[u]=V$ and $N(u) = V \setminus \{u\}$
\end{itemize}
and the symmetric differences of open neighborhoods of distinct vertices
\begin{itemize}[leftmargin=12pt, itemsep=0pt]
\item $\Delta(u,v)= (V \setminus V_i) \cup \{v\}$ for all $v \in V_i$, $1 \leq i \leq k$,
  \item $\Delta(v,v')= \{v,v'\}$ for all $v,v' \in V_i$ with $v \neq v'$ and $1 \leq i \leq k$, $n_i \geq 2$,
  \item $\Delta(v,w)= (V_i\setminus \{v\}) \cup (V_j\setminus \{w\})$ for all $v \in V_i$, $w \in V_j$ and $1 \leq i < j \leq k$.
\end{itemize}
\begin{itemize}
\item[(a)] If $n_1 = 1$ and $n_2 = \ldots = n_{\ell} = 2$ holds for some $\ell \leq k$, then
  $\Delta(v_1,w)= (\emptyset \cup (V_j\setminus \{w\}) = \{w'\}$ follows for $\{v_1\} = V_1$ and all $w \in V_j = \{w,w'\}$, $1 < j \leq \ell$. This implies $W_2 \subseteq \F^1_\X(S_k)$ for $\X \in \{\OD,\OTD\}$ and shows together with 
\begin{itemize}[leftmargin=12pt, itemsep=0pt]
\item $N[v_1]=\{v_1,u\}$ and $N(v_1) =\{u\}$, 
\item $\Delta(v,v')= \{v,v'\}$ for all $v,v' \in V_i$ with $v \neq v'$ and $\ell + 1 \leq i \leq k$
\end{itemize}
that all other hyperedges are redundant. 

    We conclude that the \OD-clutter is composed of $W_2$, $\{v_1,u\}$, and $\{v,v'\}$ for all $v,v' \in V_i$ with $v \neq v'$ and $\ell+1 \leq i \leq k$, hence $\clu_{\OD}(S_k) = K_2 \oplus 2(\ell-1)K_1 \oplus \sum_{n_i\geq 3} K_{n_i}$, whereas the \OTD-clutter is composed of $W_2$, $\{u\}$, and $\{v,v'\}$ for all $v,v' \in V_i$ with $v \neq v'$ and $\ell+1 \leq i \leq k$, and hence $\clu_{\OTD}(S_k) = (2\ell-1)K_1 \oplus \sum_{n_i\geq 3} K_{n_i}$ (note that $V^0_{\OLD}(S_k) = \{v_1\}$). 

    In both cases, this implies $\od(S_k) = \otd(S_k) = 2\ell-1 + \sum_{n_i\geq 3} (n_i-1)$.
    
  \item[(b)] If $n_1 = \ldots = n_{\ell} = 2$ for some $1 \leq \ell \leq k$, then we see that the \OD-clutter is composed of 
\begin{itemize}[leftmargin=12pt, itemsep=0pt]
\item $\Delta(v,v')= \{v,v'\}$ for all $v,v' \in V_i$ with $v \neq v'$, $1 \leq i \leq k$,
  \item $\Delta(v,w)= (V_i\setminus \{v\}) \cup (V_j\setminus \{w\})$ for all $v \in V_i$, $w \in V_j$ and $1 \leq i < j \leq \ell$,
\end{itemize}
  whereas all neighborhoods and all other symmetric differences are redundant and $V^0_{\OD}(S_k) = \{u\}$. This implies 
  $\clu_{\OD}(S_k) = K_{2\ell} \oplus \sum_{n_i\geq 3} K_{n_i}$, where $K_{2\ell}$ is composed of the vertices in $W_2$, and $\od(S_k) =  2\ell-1 + \sum_{n_i\geq 3} (n_i-1)$.

  On the other hand, the \OTD-clutter is composed of all hyperedges of $\clu_{\OD}(S_k)$ together with $N(v) = \{v',u\}$ of all vertices $v \in V_i = \{v',v\}$ and $1 \leq i \leq \ell$, whereas all other neighborhoods and symmetric differences are redundant.
  This implies $\clu_{\OTD}(S_k) = K_{2\ell+1} \oplus \sum_{n_i\geq 3} K_{n_i}$, where $K_{2\ell+1}$ is composed of the vertices in $W_2 \cup \{u\}$, and $\otd(S_k) = 2\ell + \sum_{n_i\geq 3} (n_i-1)$.
\item[(c)]  If $n_1\geq 3$, then we see that both clutters only contain 
$\Delta(v,v')= \{v,v'\}$ for all $v,v' \in V_i$ with $v \neq v'$, $1 \leq i \leq k$ as hyperedges since all neighborhoods and all other symmetric differences are redundant. Moreover, we have $V_\OD^0 (S_k) = V_\OTD^0 (S_k) = \{u\}$. This finally implies 
  $\clu_{\OD}(S_k) = \clu_{\OTD}(S_k) = \sum_{1 \leq i \leq k} K_{n_i}$ and $\od(S_k) = \otd(S_k) = \sum_{1 \leq i \leq k} (n_i-1)$ follows.  
  \end{itemize}
\end{proof}

Hence, for $k$-clique-stars $S_k$, the \OD- and the \OTD-numbers are equal if $n_1 = 1$ or if $n_1\geq 3$ for all $1 \leq i \leq k$ holds, but differ otherwise. 
In particular, calling a $k$-clique-star $S_k$ with $n_1 = \ldots = n_k = 2$ for some $k \geq 2$ to be a \emph{$k$-fan} $F_k$, the arguments from the previous proof for the case (b) show that $\clu_{\OD}(F_k)=K_{2k}$ and $V_\OD^0 (F_k) = \{u\}$ holds which implies $\od(F_k)=2k-1$, whereas $\clu_{\OTD}(F_k)=K_{2k+1}$ and $\otd(F_k)=2k$ follows. 

We note that, by Lemma \ref{lem_OD_union-cliques}(a), certain disjoint unions of cliques, including matchings $kK_2$, are examples of graphs $G$ where $\od(G)$ is \emph{larger} than the sum of the \OD-numbers of its components. 
This behavior does not apply to any other \X-number for $\X \in \{\ID, \DTD, \LD, \LTD, \OLD\}$. 

\subsection{Two families of bipartite graphs}
\begin{figure}[h]
\begin{center}
\includegraphics[scale=1.0]{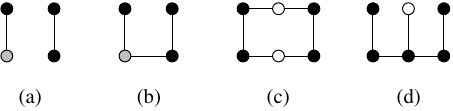}
\caption{Bipartite graphs (black vertices belong to a minimum \OD-code, grey vertices need to be added to obtain an \OTD-code), where (a) is the $2K_2$, (b) the $P_4$, (c) the $C_6$, (d) the bow.}
\label{fig_bip}
\end{center}
\end{figure}

Figure \ref{fig_bip} shows some small \otdadmis ~bipartite graphs. 
It is easy to see that $\od(G)$ and $\otd(G)$ differ for $G \in \{2K_2,P_4\}$ and are equal for $G \in \{C_6,\mbox{ bow}\}$.
We next exhibit families of bipartite graphs where the \OD- and the \OTD-numbers differ.

For any integer $k \geq 1$, the \emph{half-graph} $B_k=(U \cup W, E)$ is the bipartite graph with its stable vertex sets $U = \{u_1, \ldots , u_k\}$ and $W = \{w_1, \ldots , w_k \}$ and edges $u_iw_j$ if and only if $i \leq j$ (see Figure~\ref{fig:half_graph}). 
In particular, we have $B_1=K_2$ and $B_2=P_4$. Moreover, we clearly see that half-graphs are connected and open-twin-free and hence, are both \OD- and \otdadmis. 

In \cite{FGRS2021} it was shown that the only graphs whose \OTD-numbers equal the order of the graph are the disjoint unions of half-graphs. In particular, we have $\otd(B_k) = |V(B_k)| = 2k$. Now, let $G = (V,E)$ be a graph that is a disjoint union of half-graphs. By Theorem~\ref{thm_osd_bounds-X}(b), therefore, we have $\od(G) \geq \otd(G)-1 = |V|-1$. Moreover, by Theorem~\ref{thm_osd_bounds-n}, we have $\od(G) \leq |V|-1$. Hence, combining the two inequalities, we get the following corollary.
\begin{corollary}\label{lem_OSD_half-graphs}
For a graph $G=(V,E)$ being the disjoint union of half-graphs, we have $\od(G) = |V|-1$. In particular, for a half-graph $B_k$, we have $\od(B_k) = 2k-1$.
\end{corollary}

Corollary~\ref{lem_OSD_half-graphs} shows that, on the one hand, disjoint unions of half-graphs are further examples of graphs $G$ where $\od(G)$ is \emph{larger} than the sum of the \OD-numbers of its components and, on the other hand, that 
half-graphs and their disjoint unions are extremal examples of graphs whose \OD-numbers also attain the general upper bound in Theorem \ref{thm_osd_bounds-n}.

Moreover, we note that the upper bound from Theorem \ref{thm_osd_bounds-n} does \emph{no}t apply to \odadmis ~graphs having an isolated vertex. To see this, consider the graph $G = B_k + K_1$ for some $k \geq 1$. By Theorem \ref{thm_osd_bounds-X}(b), we have $\od(G) = \otd(B_k) + 1 = 2k+1 = |V(G)|$. 
As half-graphs and their disjoint unions are the only graphs whose \OTD-numbers equal the order of the graph by \cite{FGRS2021}, adding an isolated vertex to them yields the only graphs whose \OD-numbers equal the order of the graph.

\begin{figure}[h]
    \centering
    \begin{subfigure}{0.4\textwidth}
    \centering
    \begin{tikzpicture}[myv/.style={circle, draw, inner sep=1.7pt}, scale=1.5]
    \node[myv] (v1) at (0,0) {};
    \node[myv, fill=black] (v2) at (0.5,0) {};
    \node[myv, fill=black] (v3) at (1,0) {};
    \node[myv, fill=black] (v4) at (1.5,0) {};
    \node[myv, fill=black] (v5) at (2,0) {};
    \node[myv, fill=black] (v6) at (2.5,0) {};
    \node[myv, fill=black] (u1) at (0,0.8) {};
    \node[myv, fill=black] (u2) at (0.5,0.8) {};
    \node[myv, fill=black] (u3) at (1,0.8) {};
    \node[myv, fill=black] (u4) at (1.5,0.8) {};
    \node[myv, fill=black] (u5) at (2,0.8) {};
    \node[myv, fill=black] (u6) at (2.5,0.8) {};
    \draw (v1) -- (u1);
    \draw (v1) -- (u2);
    \draw (v1) -- (u3);
    \draw (v1) -- (u4);
    \draw (v1) -- (u5);
    \draw (v1) -- (u6);
    \draw (v2) -- (u2);
    \draw (v2) -- (u3);
    \draw (v2) -- (u4);
    \draw (v2) -- (u5);
    \draw (v2) -- (u6);
    \draw (v3) -- (u3);
    \draw (v3) -- (u4);
    \draw (v3) -- (u5);
    \draw (v3) -- (u6);
    \draw (v4) -- (u4);
    \draw (v4) -- (u5);
    \draw (v4) -- (u6);
    \draw (v5) -- (u5);
    \draw (v5) -- (u6);
    \draw (v6) -- (u6);
\end{tikzpicture}
    \caption{Half-graph $B_6$.}
    \label{fig:half_graph}
    \end{subfigure}
    \hspace{2mm}
    \begin{subfigure}{0.4\textwidth}
    \centering
    \begin{tikzpicture}[myv/.style={circle, draw, inner sep=1.7pt}, scale=1.5]
    \node[myv] (u0) at (1.25,1) {};
    \node[myv, fill=black] (u1) at (0,0.5) {};
    \node[myv, fill=black] (u2) at (0.5,0.5) {};
    \node[myv, fill=black] (u3) at (1,0.5) {};
    \node[myv, fill=black] (u4) at (1.5,0.5) {};
    \node[myv, fill=black] (u5) at (2,0.5) {};
    \node[myv, fill=black] (u6) at (2.5,0.5) {};
    \node[myv] (w1) at (0,0) {};
    \node[myv, fill=black] (w2) at (0.5,0) {};
    \node[myv, fill=black] (w3) at (1,0) {};
    \node[myv, fill=black] (w4) at (1.5,0) {};
    \node[myv, fill=black] (w5) at (2,0) {};
    \node[myv, fill=black] (w6) at (2.5,0) {};
    \draw (u0) -- (u1);
    \draw (u0) -- (u2);
    \draw (u0) -- (u3);
    \draw (u0) -- (u4);
    \draw (u0) -- (u5);
    \draw (u0) -- (u6);
    \draw (w1) -- (u1);
    \draw (w2) -- (u2);
    \draw (w3) -- (u3);
    \draw (w4) -- (u4);
    \draw (w5) -- (u5);
    \draw (w6) -- (u6);
    \end{tikzpicture}
    \caption{6-double star $D_6$.}
    \label{fig:double_star}
    \end{subfigure}
    \caption{The black vertices depict an \OD-code of the respective graph.}
\end{figure}

To give an independent proof for the \OD- and \OTD-numbers of half-graphs, we next provide their \X-clutters.

\begin{lemma} \label{lem_OD-OTD_half-graphs}
For any half-graph $B_k$ with $k \geq 2$, we have
\begin{itemize}[leftmargin=12pt, itemsep=0pt]
\item $\clu_{\OD}(B_k) = (2k-2)K_1 \oplus K_{2}$ and $\od(B_k) =  2k-1$ but
\item $\clu_{\OTD}(B_k) = 2kK_1$ and $\otd(B_k) = 2k$.
\end{itemize}
\end{lemma}

\begin{proof}
  Consider a half-graph $B_k=(U \cup W, E)$ with its stable vertex sets $U = \{u_1, \ldots , u_k\}$ and $W = \{w_1, \ldots , w_k \}$ and edges $u_iw_j$ if and only if $i \leq j$ for some $k \geq 2$.
  In order to determine $\clu_{\X}(B_k)$ and $\x(B_k)$ for $\X \in \{\OD,\OTD\}$, we first construct the \X-hypergraphs. Recall that the following hyperedges are involved: the closed or open neighborhoods
\begin{itemize}[leftmargin=12pt, itemsep=0pt]
\item $N[u_i] = \{u_i\} \cup \{w_i, \ldots, w_k\}$ and $N(u_i) = \{w_i, \ldots, w_k\}$ for $1 \leq i \leq k$,
\item $N[w_j] = \{u_1, \ldots , u_j\} \cup \{w_j\}$ and $N(w_j) = \{u_1, \ldots , u_j\}$ for $1 \leq j \leq k$,
\end{itemize}
and the symmetric differences of open neighborhoods of distinct vertices
\begin{itemize}[leftmargin=12pt, itemsep=0pt]
\item $\Delta(u_i,u_j) = \{w_i, \ldots, w_{j-1}\}$ for $1 \leq i < j \leq k$,
\item $\Delta(u_i,w_j) = \{w_i, \ldots, w_k\} \cup \{u_1, \ldots , u_j\}$ for $1 \leq i,j \leq k$,
\item $\Delta(w_i,w_j) = \{u_{1+1}, \ldots , u_j\}$ for $1 \leq i < j \leq k$.
\end{itemize}
In particular, we see that
\begin{itemize}[leftmargin=12pt, itemsep=0pt]
\item $\Delta(u_i,u_{1+1}) = \{w_i\}$ for $1 \leq i < k$,
\item $\Delta(w_i,w_{1+1}) = \{u_{1+1}\}$ for $1 \leq i < k$
\end{itemize}
belong to $\for_{\X}(B_k)$ for $\X \in \{\OD,\OTD\}$.
Noting further that $N(u_k) = \{w_k\}$ and $N(w_1) = \{u_1\}$ holds implies $\for_{\OTD}(B_k) = U \cup W$ and thus $\clu_{\OTD}(B_k) = 2kK_1$ and $\otd(B_k) = 2k$.
On the other hand, it is easy to see that $\Delta(u_k,w_1) = \{w_k,u_1\}$ is the only hyperedge in $\fac_{\OD}(B_k)$, thus $\clu_{\OD}(B_k) = (2k-2)K_1 \oplus K_{2}$ and $\od(B_k) =  2k-1$ follows.
\end{proof}

A \emph{$k$-double star} $D_k=(U \cup W, E)$ is the bipartite graph with its stable vertex sets $U = \{u_0,u_1, \ldots , u_k\}$ and $W = \{w_1, \ldots , w_k \}$ and edges $u_iw_i$ and $u_0w_i$ for all $w_i \in W$ (see Figure~\ref{fig:double_star}). 
Then, we have $D_1=P_3$ and $D_2=P_5$. Moreover, we clearly see that $k$-double stars with $k \geq 2$ are connected and open-twin-free and hence, are both \OD- and \otdadmis. As the next lemma shows, $k$-double stars also provide examples of bipartite graphs where the \OD- and the \OTD-numbers differ.
\begin{lemma} \label{lem_OSD_double-stars} For a $k$-double star $D_k$ with $k \geq 3$,
  $\clu_{\OD}(D_k)$ is composed of
\begin{itemize}[leftmargin=12pt, itemsep=0pt]
  \item $\{u_i,w_i\}$ for $1 \leq i \leq k$,
  \item $\{u_0,u_j,w_i\}$ for $1 \leq i,j \leq k$ for $i \neq j$,
  \item $\{u_i,u_j\}$ and $\{w_i,w_j\}$ for $1 \leq i < j \leq k$,
\end{itemize}  
and thus $\od(D_k) = 2k-1$ follows, whereas $\clu_{\OTD}(D_k)$ is composed of
\begin{itemize}[leftmargin=12pt, itemsep=0pt]
  \item $\{w_i\}$ and $\{u_0,u_i\}$ for $1 \leq i < k$,
  \item $\{u_i,u_j\}$ for $1 \leq i < j \leq k$.
\end{itemize} 
and $\otd(D_k) = 2k$ holds.
\end{lemma}

\begin{proof}
Consider a $k$-double star $D_k=(U \cup W, E)$ with vertex sets $U = \{u_0,u_1, \ldots , u_k\}$ and $W = \{w_1, \ldots , w_k \}$ and edges $u_i w_i$ and $u_0 w_i$ for all $w_i \in W$ for some $k \geq 3$. 
In order to determine $\mathcal{C}_{\X}(D_k)$ and $\x(D_k)$ for $\X \in \{\OD,\OTD\}$, we first construct the \X-hypergraphs. Recall that the following hyperedges are involved: the closed or open neighborhoods
\begin{itemize}[leftmargin=12pt, itemsep=0pt]
\item $N[u_0] = \{u_0\} \cup W$ and $N(u_0) = W$,
\item $N[u_i] = \{u_i,w_i\}$ and $N(u_i) = \{w_i\}$ for $1 \leq i \leq k$,
\item $N[w_i] = \{u_0,u_i,w_i\}$ and $N(w_i) = \{u_0,u_i\}$ for $1 \leq i \leq k$,
\end{itemize}
and the symmetric differences of open neighborhoods of distinct vertices
\begin{itemize}[leftmargin=12pt, itemsep=0pt]
\item $\Delta(u_0,u_i) = W \setminus \{w_i\}$ for $1 \leq i \leq k$,
\item $\Delta(u_0,w_i) = W \cup \{u_0,u_i\}$ for $1 \leq i \leq k$,
\item $\Delta(u_i,u_j) = \{w_i,w_j\}$ for $1 \leq i < j \leq k$,
\item $\Delta(u_i,w_j) = \{u_0,u_j,w_i\}$ for $1 \leq i,j \leq k$,
\item $\Delta(w_i,w_j) = \{u_i,u_j\}$ for $1 \leq i < j \leq k$.
\end{itemize}
The \OD-hypergraph involves the closed neighborhoods and all symmetric differences, including
\begin{itemize}[leftmargin=12pt, itemsep=0pt]
\item $N[u_i] = \{u_i,w_i\}$ for $1 \leq i \leq k$,
\item $\Delta(u_i,u_j) = \{w_i,w_j\}$ for $1 \leq i < j \leq k$,
\item $\Delta(u_i,w_j) = \{u_0,u_j,w_i\}$ for $1 \leq i,j \leq k$ for $i \neq j$,
\item $\Delta(w_i,w_j) = \{u_i,u_j\}$ for $1 \leq i < j \leq k$,
\end{itemize}
which implies that all other hyperedges are redundant for any $k \geq 3$.
Analyzing the hyperedges remaining in $\clu_{\OD}(D_k)$ shows that $(U \setminus \{u_i\}) \cup (W \setminus \{w_j\})$ for $0 \leq i \leq k$, $1 \leq j \leq k$ and $i \neq j$ as well as $(U \setminus \{u_0,u_i\}) \cup W$ for $1 \leq i \leq k$ are the only minimum \OD-codes of $D_k$ and thus $\od(D_k) = 2k-1$ follows.

The \OTD-hypergraph involves the open neighborhoods and all symmetric differences, including
\begin{itemize}[leftmargin=12pt, itemsep=0pt]
\item $N(u_i) = \{w_i\}$ for $1 \leq i \leq k$,
\item $N(w_i) = \{u_0,u_i\}$ for $1 \leq i \leq k$,
\item $\Delta(w_i,w_j) = \{u_i,u_j\}$ for $1 \leq i < j \leq k$,
\end{itemize}
which implies that all other hyperedges are redundant for any $k \geq 3$.
Analyzing the hyperedges remaining in $\clu_{\OTD}(D_k) = kK_1 \oplus K_{k+1}$ shows that $(U \setminus \{u_i\}) \cup W$ for $0 \leq i \leq k$ are the only minimum \OTD-codes of $D_k$ and $\otd(D_k) = 2k$ holds.
\end{proof}
  
\subsection{Families of split graphs}
\begin{figure}[h]
\begin{center}
\includegraphics[scale=1.0]{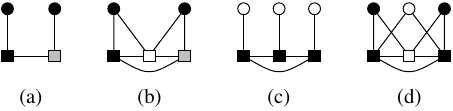}
\caption{Split graphs (the squares indicate the vertices in $Q$, black vertices belong to a minimum \OD-code, grey vertices need to be added to obtain an \OTD-code), where (a) is the $P_4$, (b) the gem, (c) the net, (d) the sun.}
\label{fig_split}
\end{center}
\end{figure}

In order to study \OD-codes of split graphs and compare them with the \OTD-codes, we restrict ourselves to split graphs $G$ without open twins and isolated vertices. 
This further implies that $G$ is connected and $Q$ non-empty (as, otherwise, every component not containing the clique $Q$ needs to be an isolated vertex from $S$, contradicting our assumptions).
Figure \ref{fig_split} shows some small \otdadmis ~split graphs. 
It is easy to see that $\od(G)$ and $\otd(G)$ differ for $G \in \{P_4, \mbox{ gem}\}$ and are equal for $G \in \{\mbox{net},\mbox{ sun}\}$.

We next examine \OD-codes in two families of split graphs for which the exact \OTD-numbers are known from \cite{ABLW_2022}.
A \emph{headless spider} is a split graph with $Q = \{q_1, \ldots, q_k\}$ and $S = \{s_1, \ldots, s_k\}$. In addition, a headless spider is \emph{thin} (respectively, \emph{thick}) if $s_i$ is adjacent to $q_j$ if and only if $i = j$ (respectively, $i \neq j$). By definition, it is clear that the complement of a thin headless spider $H_k$ is a thick headless spider $\overline H_k$, and vice-versa. 
We have $H_2 = \overline H_2 = P_4$, the two headless spiders $H_3 = net$ and $\overline H_3 = sun$ are depicted in Figures~\ref{fig_split}(c) and \ref{fig_split}(d), respectively. Moreover, it is easy to check that thin and thick headless spiders have no twins.

In \cite{ABLW_2022}, it was shown that $\otd(H_k) = k$ for $k \geq 3$ and $\otd(\overline H_k) = k+1$ for $k \geq 3$. We next analyse the \OD-numbers of thin and thick headless spiders.

\begin{lemma}\label{lem_ThickSpider_OD_clutter}
For a thick headless spider $\overline H_k=(Q \cup S, E)$ with $k \geq 3$, the \OD-clutter is composed of the symmetric differences of non-adjacent vertices
\begin{itemize}[leftmargin=12pt, itemsep=0pt]
  \item $N(q_i) \Delta N(s_i) = S \setminus \{s_i\}$ for all $i \in \{1, \ldots, k\}$,
  \item $N(s_i) \Delta N(s_j) = \{q_i,q_j\}$ for $i \in \{1, \ldots, k\}$, $j \in \{1, \ldots, k\} \setminus \{i\}$
\end{itemize}
and we thus have $\clu(\overline H_k) =  {\cal{R}}_{|S|}^{|S|-1} \cup {\cal{R}}_{|Q|}^2$ and $\od(\overline H_k) = k+1$.
\end{lemma}

\begin{proof}
Consider a thick headless spider $\overline H_k=(Q \cup S, E)$ with $k \geq 3$. 
$\hyp_{\OD}(\overline H_k)$ is composed of the closed neighborhoods
\begin{itemize}[leftmargin=12pt, itemsep=0pt]
\item $N[s_i] = Q \setminus \{q_i\} \cup \{s_i\}$ for all $s_i \in S$,
\item $N[q_i] = Q \cup \{s_i\}$ for all $q_i \in Q$
\end{itemize}
and the symmetric differences
\begin{itemize}[leftmargin=12pt, itemsep=0pt]
\item $\Delta(s_i,s_j) = \{q_i,q_j\}$ for all distinct $s_i,s_j \in S$, 
\item $\Delta(q_i,q_j) = \{q_i,q_j\} \cup \{s_i,s_j\}$ for all distinct $q_i,q_j \in C$,
\item $\Delta(q_i,s_j) =
\begin{cases}
  S \setminus \{s_i\}, &\mbox{if } i=j; \\
  \{q_i,q_j\} \cup S \setminus \{s_i\}, &\mbox{if } i \neq j.\\
\end{cases}$
\end{itemize}
This shows that all neighborhoods are redundant as well as all $\Delta(q_i,q_j)$ and $\Delta(q_i,s_j)$ with $i \neq j$ so that only $\Delta(s_i,s_j)$ and $\Delta(q_i,s_i)$ belong to $\clu_{\OD}(\overline H_k)$.
Hence, we obtained $\mathcal{C}_{\OD}(\overline H_k) = {\cal{R}}_{|S|}^{|S|-1} \cup {\cal{R}}_{|Q|}^2$.
Finally, this implies $\od(\overline H_k) = 2 + (k-1) = k+1$ by the result from \cite{S_1989}.
\end{proof}

\begin{lemma}\label{lem_ThinSpider_OD_clutter}
For a thin headless spider $H_k=(Q \cup S, E)$ with $k \geq 4$, the \OD-clutter is composed of 
\begin{itemize}[leftmargin=12pt, itemsep=0pt]
  \item $N[s_i] = \{q_i,s_i\}$ for all $i \in \{1, \ldots, k\}$,
  \item $N(s_i) \Delta N(s_j) = \{q_i,q_j\}$ for $i \in \{1, \ldots, k\}$, $j \in \{1, \ldots, k\} \setminus \{i\}$
\end{itemize}
and we thus have $\mathcal{C}_{\OD}(H_k) = H_k$ and $\od(H_k) = k$.
\end{lemma}

\begin{proof}
Consider a thin headless spider $H_k=(Q \cup S, E)$ with $k \geq 4$. 
$\hyp_{\OD}(H_k)$ is composed of the closed neighborhoods
\begin{itemize}[leftmargin=12pt, itemsep=0pt]
\item $N[s_i] = \{s_i,q_i\}$ for all $s_i \in S$,
\item $N[q_i] = Q \cup \{s_i\}$ for all $q_i \in Q$
\end{itemize}
and the symmetric differences
\begin{itemize}[leftmargin=12pt, itemsep=0pt]
\item $\Delta(s_i,s_j) = \{q_i,q_j\}$ for all distinct $s_i,s_j \in S$, 
\item $\Delta(q_i,q_j) = \{q_i,q_j\} \cup \{s_i,s_j\}$ for all distinct $q_i,q_j \in Q$,
\item $\Delta(q_i,s_j) =
\begin{cases}
  Q \cup \{s_i\}, &\mbox{if } i=j; \\
  (Q \setminus \{q_i,q_j\}) \cup \{s_i\}, &\mbox{if } i \neq j.\\
\end{cases}$
\end{itemize}
We deduce that in $\clu_{\OD}(H_k)$ with $k \geq 4$ only
\begin{itemize}[leftmargin=12pt, itemsep=0pt]
\item $N[s_i] = \{s_i,q_i\}$ for all $s_i \in S$,
\item $\Delta(s_i,s_j) = \{q_i,q_j\}$ for all distinct $s_i,s_j \in S$
\end{itemize}
remain and thus $\mathcal{C}_{\OD}(H_k) = H_k$ follows.
Hence $\mathcal{C}_{\OD}(H_k)$ contains the disjoint hyperedges $\{q_i,s_i\}$ for $1 \leq i \leq k$ which shows $\od(H_k) \geq k$.
On the other hand, $Q$ is clearly a cover of $\mathcal{C}_{\OD}(H_k) = H_k$ of cardinality $k$, which finally implies $\od(H_k) = k$.
\end{proof}

Moreover, it is easy to see that $\od(H_3) = 3$ (see Figure \ref{fig_split}(c)), hence the two lemmas combined with the results from \cite{ABLW_2022} show that for thin and thick headless spiders $H_k$ and $\overline H_k$, respectively, the \OD- and the \OTD-numbers are equal for all $k \geq 3$.

In \cite{CW_ISCO2024}, the question was raised whether there also exist infinite families of open twin-free split graphs where the \OD- and the \OTD-numbers differ.
To answer this question, we consider the following extension of thin headless spiders.
Consider a split graph $G=(Q \cup S, E)$ with $Q = \{q_1, \ldots, q_k\}$ and $S = \{s_0,s_1, \ldots, s_k\}$,
where $s_i$ is adjacent to $q_i$ for all $1 \leq i \leq k$ and $s_0$ is adjacent to $\{q_1, \ldots, q_{k-1}\}$.
By definition, it is clear that $G[Q \cup S\setminus\{s_0\}]$
is isomorphic to the thin headless spider $H_k$. 
Hence, we call such split graphs \emph{extended thin spiders} $E_k$.

\begin{figure}[h]
\begin{center}
\includegraphics[scale=1.0]{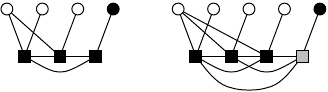}
\caption{Extended thin spiders $E_3$ and $E_4$ (the squares indicate the vertices in $Q$, black vertices belong to a minimum \OD-code, the grey vertex needs to be added to obtain an \OTD-code).}
\label{fig_ext-thin-spiders}
\end{center}
\end{figure}

It is easy to see that for $k \geq 3$, extended thin spiders $E_k$ have neither isolated vertices nor open twins, see Figure \ref{fig_ext-thin-spiders} for illustration. 
We can show:

\begin{lemma}\label{lem_ExtThinSpider_OD-OTD_clutter}
For an extended thin spider $E_k=(Q \cup S, E)$ with $k \geq 4$, the 
\begin{itemize}[leftmargin=12pt, itemsep=0pt]
\item \OD-clutter is composed of $\{s_k\}$, $\{q_i,s_i\}$ for all $1 \leq i < k$, and $\{q_i,q_j\}$ for all $1 \leq i < j \leq k$ and we thus have $\od(E_k) = k$;
\item \OTD-clutter is composed of $\{s_k\}$ and $\{q_i\}$ for $1 \leq i \leq k$ and we have $\otd(E_k) = k+1$.
\end{itemize}
\end{lemma}

\begin{proof}
  Consider an extended thin spider $E_k=(Q \cup S, E)$ with $k \geq 4$.
  In order to determine $\mathcal{C}_{\X}(E_k)$ and $\x(E_k)$ for $\X \in \{\OD,\OTD\}$, we first construct the \X-hypergraphs. Recall that the following hyperedges are involved: the closed or open neighborhoods
\begin{itemize}[leftmargin=12pt, itemsep=0pt]
  \item $N[s_0] = \{s_0\} \cup (Q \setminus\{q_k\})$ and $N(s_0) =Q \setminus\{q_k\}$,
\item $N[s_i] = \{s_i,q_i\}$ and $N(s_i) =\{q_i\}$ for $1 \leq i \leq k$,
\item $N[q_i] = Q \cup \{s_0,s_i\}$ and $N(q_i) =(Q \setminus\{q_i\}) \cup \{s_0,s_i\}$ for $1 \leq i < k$,
\item $N[q_k] = Q \cup \{s_k\}$ and $N(q_k) =(Q \setminus\{q_k\}) \cup \{s_k\}$,
\end{itemize}
and the symmetric differences of open neighborhoods of distinct vertices
\begin{itemize}[leftmargin=12pt, itemsep=0pt]
\item $\Delta(s_0,s_i) = Q \setminus\{q_i,q_k\}$ for $1 \leq i < k$,
\item $\Delta(s_0,s_k) = Q$,
\item $\Delta(s_0,q_i) = \{q_i,q_k\} \cup \{s_0,s_i\}$ for $1 \leq i < k$,
\item $\Delta(s_0,q_k) = \{s_k\}$,
  
\item $\Delta(s_i,s_j) = \{q_i,q_j\}$ for $1 \leq i < j \leq k$,
\item $\Delta(s_i,q_i) = Q \cup \{s_0,s_i\}$ for $1 \leq i < k$,
\item $\Delta(s_k,q_k) = Q \cup \{s_k\}$,
\item $\Delta(s_i,q_j) = (Q \setminus \{q_i,q_j\}) \cup \{s_0,s_j\}$ for $1 \leq j < i \leq k$,
\item $\Delta(s_i,q_k) = (Q \setminus \{q_i,q_k\}) \cup \{s_k\}$ for $1 \leq i < k$,
  
\item $\Delta(q_i,q_j) = \{q_i,q_j\} \cup \{s_i,s_j\}$ for $1 \leq j < i \leq k$,
\item $\Delta(q_i,q_k) = \{q_i,q_k\} \cup \{s_0,s_i,s_k\}$ for $1 \leq i < k$.
\end{itemize}
The \OD-hypergraph involves the closed neighborhoods and all symmetric differences, including
\begin{itemize}[leftmargin=12pt, itemsep=0pt]
\item $\Delta(s_0,q_k) = \{s_k\}$,
\item $N[s_i] = \{s_i,q_i\}$ for $1 \leq i < k$, and 
\item $\Delta(s_i,s_j) = \{q_i,q_j\}$ for $1 \leq i < j \leq k$,
\end{itemize}
which implies that all other hyperedges are redundant and $V^0_{\OD}(E_k)=\{s_0\}$. The subset $\{s_k\}$ and $\{s_i,q_i\}$ for $1 \leq i < k$ of disjoint hyperedges shows $\od(E_k) \geq k$. On the other hand, $Q \setminus \{q_k\} \cup \{s_k\}$ is clearly a cover of $\mathcal{C}_{\OD}(E_k)$ of cardinality $k$, which finally implies $\od(E_k) = k$.

The \OTD-hypergraph involves the open neighborhoods and all symmetric differences, including
\begin{itemize}[leftmargin=12pt, itemsep=0pt]
\item $N(s_i) =\{q_i\}$ for $1 \leq i \leq k$,
\item $\Delta(s_0,q_k) = \{s_k\}$,
\end{itemize}
which implies that all other hyperedges are redundant and $V^0_{\OLD}(E_k)=S \setminus \{s_k\}$. This clearly shows that $Q \cup \{s_k\}$ is the only minimum cover of $\mathcal{C}_{\OTD}(E_k)$ and we have $\otd(E_k) = k+1$.
\end{proof}

\subsection{Families of thin suns}
\begin{figure}[!t]
\begin{center}
\includegraphics[scale=1.0]{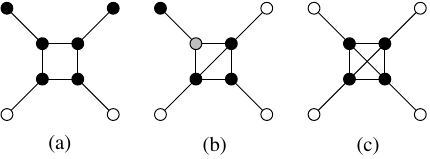}
\caption{The three thin suns $T_4$ where (a) is a sunlet and (c) a thin headless spider (black vertices belong to a minimum \OD-code, the grey vertex needs to be added to obtain an \OTD-code).}
\label{Fig_thin_suns}
\end{center}
\end{figure}

The result on thin headless spiders can be further generalized to thin suns. 
A \emph{sun} is a graph $G=(C\cup S,E)$ whose vertex set can be partitioned into $S$ and $C$, where, for an integer $k \geq 3$, the set $S = \{s_1, \ldots, s_k\}$ is a stable set and $C = \{c_1, \ldots, c_k\}$ is a (not necessarily chordless) cycle. 
A \textit{thin sun} $T_k=(C\cup S,E)$ is a sun where $s_i$ is adjacent to $c_j$ if and only if $i = j$. Therefore, thin headless spiders are special thin suns where all chords of the cycle $C$ are present (such that $C$ induces a clique). 
Other special cases of thin suns are \emph{sunlets} where no chords of the cycle $C$ are present (such that $C$ induces a hole). 
For illustration, for $k=3$, the (only) thin sun $T_3$ equals the thin headless spider $H_3$ (see Figure~\ref{fig_split}(c)); for $k=4$, the three possible thin suns $T_4$ are depicted in Figure~\ref{Fig_thin_suns}.
We further note that thin suns have neither isolated vertices nor open twins. 
We first determine the \OD-clutter of thin suns.

\begin{lemma}\label{lem_ThinSun_clutter} 
For a thin sun $T_k=(C \cup S, E)$ with $k \geq 4$, the \OD-clutter $\clu_{\OD}(T_k)$ is composed of 
\begin{itemize}[leftmargin=12pt, itemsep=0pt]
\item $N[s_i] = \{s_i,c_i\}$ for all $s_i \in S$,
\item $\Delta(s_i,s_j) = \{c_i,c_j\}$ for all distinct $s_i,s_j \in S$, 

\item $\Delta(c_i,c_j) = 
\begin{cases}
  \{s_i,s_j\}, &\mbox{if } c_i, c_j \mbox{ are open $C$-twins}; \\
  \{s_i,s_j,c_{\ell}\}, &\mbox{if } c_i, c_j \mbox{ non-adjacent and } \{c_{\ell}\} = N_C(c_i) \Delta N_C(c_j),
\end{cases}$

\item $\Delta(c_i,s_j) = \{s_i,c_{\ell}\}$ if $c_i, c_j$ are adjacent and $\{c_{\ell}\} = N_C(c_i) \setminus \{c_{j}\}$, $\ell \neq i,j$.
\end{itemize}
\end{lemma}

\begin{proof}
$\hyp_{\OD}(T_k)$ is composed of the closed neighborhoods
\begin{itemize}[leftmargin=12pt, itemsep=0pt]
\item $N[s_i] = \{s_i,c_i\}$ for all $s_i \in S$,
\item $N[c_i] = N_C[c_i] \cup \{s_i\}$ for all $c_i \in C$
\end{itemize}
and the symmetric differences
\begin{itemize}[leftmargin=12pt, itemsep=0pt]
\item $\Delta(s_i,s_j) = \{c_i,c_j\}$ for all distinct $s_i,s_j \in S$, 
\item $\Delta(c_i,c_j) = (N_C(c_i) \Delta N_C(c_j)) \cup \{s_i,s_j\}$ for all distinct $c_i,c_j \in C$,
\item $\Delta(c_i,s_j) =
\begin{cases}
  N[c_i], &\mbox{if } i=j; \\
  (N_C(c_i) \cup \{s_i\}) \Delta \{c_j\}, &\mbox{if } i \neq j.\\
\end{cases}$
\end{itemize}
This shows that $N[s_i]$ and $\Delta(s_i,s_j)$ belong to $\clu_{\OD}(T_k)$ which further implies that the following hyperedges from $\hyp_{OD}(T_k)$ are redundant:
\begin{itemize}[leftmargin=12pt, itemsep=0pt]
\item $N[c_i] = N_C[c_i] \cup \{s_i\}$ for all $c_i \in C$,
\item $\Delta(c_i,c_j)$ if $|N_C(c_i) \Delta N_C(c_j)| \geq 2$,
\item $\Delta(c_i,s_j)$ if $c_i, c_j$ are non-adjacent or if $c_i, c_j$ are adjacent but $|N_C(c_i)| \geq 3$.
\end{itemize}
In the remaining cases, $\Delta(c_i,c_j)$ and $\Delta(c_i,s_j)$ belong to $\clu_{\OD}(T_k)$.
\end{proof}

We call two vertices $c_i$ and $c_j$ of a thin sun $T_k=(C\cup S,E)$ \emph{open $C$-twins} if $c_i$ and $c_j$ are non-adjacent and $N_C(c_i) = N_C(c_j)$, where $N_C(v) = N(v) \cap C$. 
For instance, the sunlet in Figure~\ref{Fig_thin_suns}(a) and the thin sun in Figure~\ref{Fig_thin_suns}(b) have open $C$-twins, whereas the thin headless spider in Figure~\ref{Fig_thin_suns}(c) does not.

In \cite{ABLW_2022}, it was shown that for a thin sun $T_k$ with $k \geq 4$ and without open $C$-twins, the set $C$ is the unique minimum \OTD-code of $T_k$ and thus, we have $\otd(T_k) = k$. Now, with regards to \OD-numbers of thin suns, we show the following.

\begin{lemma}\label{lem_ThinSun_1} 
For a thin sun $T_k=(C \cup S, E)$ with $k \geq 4$ and without open $C$-twins, 
$C$ is a minimum $\OD$-code and hence, we have $\od(T_k)=|C|=k$.
\end{lemma}

\begin{proof}
Consider a thin sun $T_k=(C \cup S, E)$ with $k \geq 4$ and without open $C$-twins. 
By Lemma \ref{lem_ThinSun_clutter}, the \OD-clutter $\clu_{\OD}(T_k)$ contains $\{s_i,c_i\}$ for all $s_i \in S$, which implies the lower bound $\od(T_k) \geq k$. 
On the other hand, it is easy to verify that, by Lemma \ref{lem_ThinSun_clutter}, all hyperedges of the \OD-clutter $\clu_{\OD}(T_k)$ have a nonempty intersection with $C$ (when $T_k$ has no open $C$-twins) so that $C$ is a cover of $\clu_{\OD}(T_k)$ of cardinality $k$. 
Hence, $C$ is a minimum \OD-code and the assertion $\gamma^{\OD}(T_k)=|C|=k$ follows.
\end{proof}

Therefore, thin suns without open $C$-twins are examples of graphs where the \OD- and the \OTD-number are equal. 
This applies in particular to sunlets and to thin headless spiders. However, for thin suns $T_k$ with open $C$-twins, $\od(T_k)$ and $\otd(T_k)$ may differ. For instance, for the thin sun $T_4$ depicted in Fig. \ref{Fig_thin_suns}(b), it can be checked that $\od(T_4) = 4 < 5 = \otd(T_4)$. We call a thin sun $T_k=(C \cup S, E)$ \emph{almost complete} if $k = 2\ell$ and $c_i$ is non-adjacent to $c_{i + \ell}$ but is adjacent to all other $c_j \in C$.
For instance, the sunlet in Fig. \ref{Fig_thin_suns}(a) is almost complete. 
We can show:

\begin{lemma}\label{lem_ThinSun_2}
For an almost complete thin sun $T_{2\ell}$ with $\ell \geq 3$, we have $\od(T_{2\ell}) = 3\ell - 1$ and $\otd(T_{2\ell}) = 3\ell$. 
\end{lemma}

\begin{proof}
Consider an almost complete thin sun $T_{2\ell}$ with $\ell \geq 3$ where $c_i,c_{i+\ell}$ form open $C$-twins. 
By Lemma \ref{lem_ThinSun_clutter}, the \OD-clutter $\clu_{\OD}(T_{2\ell})$ is composed of 
\begin{itemize}[leftmargin=12pt, itemsep=0pt]
\item $N[s_i] = \{s_i,c_i\}$ for all $s_i \in S$,
\item $\Delta(s_i,s_j) = \{c_i,c_j\}$ for all distinct $s_i,s_j \in S$, 
\item $\Delta(c_i,c_{i+\ell}) = \{s_i,s_{i+\ell}\}$ for all $c_i \in C$.
\end{itemize}
Hence, every cover of $\clu_{\OD}(T_{2\ell})$ has to contain at least all but one vertex from $C$ and one of $\{s_i,s_{i+\ell}\}$ for all $1 \leq i \leq \ell$, showing that $\od(T_{2\ell}) \geq 3\ell - 1$. 
We next observe that $V' = \{s_1, \ldots, s_{\ell}\} \cup C \setminus \{c_j\}$ for some $j \in \{1, \ldots, {\ell}\}$ is a cover of $\clu_{\OD}(T_{2\ell})$. Indeed, $V'$ meets all
\begin{itemize}[leftmargin=12pt, itemsep=0pt]
\item $N[s_i] = \{s_i,c_i\}$ in $s_i$ for $1 \leq i \leq \ell$ and in $c_i$ for $\ell+1 \leq i \leq 2\ell$,
\item $\Delta(s_i,s_j) = \{c_i,c_j\}$ as $V'$ contains all but one vertex from $C$, 
\item $\Delta(c_i,c_{i+\ell}) = \{s_i,s_{i+\ell}\}$ in $s_i$ for $1 \leq i \leq \ell$.
\end{itemize}
Hence, $V'$ is an \OD-code of $T_{2\ell}$ and $\od(T_{2\ell}) = 3\ell - 1$ follows. 
Furthermore, from \cite{ABLW_2022} we deduce that the \OTD-clutter $\clu_{\OTD}(T_{2\ell})$ is composed of 
\begin{itemize}[leftmargin=12pt, itemsep=0pt]
\item $N(s_i) = \{c_i\}$ for all $s_i \in S$,
\item $\Delta(c_i,c_{i+\ell}) = \{s_i,s_{i+\ell}\}$ for all $c_i \in C$.
\end{itemize}
For an almost complete thin sun $T_{2\ell}$ with $\ell \geq 3$, this implies that $C \cup \{s_1, \ldots, s_{\ell}\}$ is a minimum \OTD-code of $T_{2\ell}$ and, hence, $\otd(T_{2\ell}) = 3\ell$ follows.
\end{proof}

Hence, there exist infinitely many thin suns with open $C$-twins for which the \OD- and the \OTD-numbers differ.

\section{Polyhedra associated with \OD-codes}
\label{sec5}

As polyhedral methods turned out to be successful for many NP-hard combinatorial optimization problems in the literature, it was suggested in e.g. \cite{ABLW_2018,ABLW_2022} to apply such techniques to identification problems. 
For that, the reformulation of the studied \X-problem in a graph $G$ as covering problem in the hypergraph $\hyp_X(G)$ was used. In fact, the incidence matrix of $\hyp_X(G)$ defines the constraint system of the resulting covering problem.
We next provide some preliminaries and then study \OD-codes in this context.
 
\paragraph{About covering polyhedra.}
Consider a $0/1$-matrix $M$ with $n$ columns.
A \emph{cover} of $M$ is a vector $\mathbf x \in \mathbf Z_+^n$ such that $M \mathbf x \geq \mathbf 1$ holds where $\mathbf 1$ is the vector having $1$-entries only (note that a minimum cover $\mathbf x$ has $0/1$-entries only). 
The \emph{covering number} of $M$ can be calculated by \[
\begin{array}{rcl}
\displaystyle \tau(M) = \min \mathbf 1^T\mathbf x \ &   &    \\
M \ \mathbf x \ & \geq & \mathbf 1   \\
\mathbf x \ & \in  & \mathbf Z_+^n.  \\
\end{array}
\]
The \emph{covering polyhedron} of $M$ is defined by
$P(M)=\mbox{conv}\left\{\mathbf{x} \in \mathbf{Z}_+^n: M \mathbf x \geq \mathbf 1 \right\}$ 
and encodes the convex hull of all covers of $M$. 
The following properties of general covering polyhedra are known from \cite{BN_1989}:

\begin{theorem}[\cite{BN_1989}] \label{thm_BN}
Consider a $0/1$-matrix $M$ with $n$ columns and the covering polyhedron $P(M)$. 
\begin{itemize}[leftmargin=12pt, itemsep=0pt]
\item $P(M)$ has dimension $n$ if and only if $M$ has at least two ones per row.
\item The only facet-defining inequalities of $P(M)$ with integer coefficients and right hand side equal to one are those of the constraint system $M \mathbf x \geq \mathbf 1$.
\end{itemize}
\end{theorem}
The linear relaxation
$Q(M)=\left\{\mathbf x \in \mathbf{R}_+^n: M \mathbf x \geq \mathbf 1 \right\}$
of $P(M)$ satisfies $P(M) \subseteq Q(M)$ in general, and $M$ is called \emph{ideal} if $P(M) = Q(M)$ holds.
For all non-ideal matrices $M$, further constraints (with right hand side $> 1$) have to be added to the system $M \mathbf x \geq \mathbf 1$ in order to describe $P(M)$ using real variables instead of integral ones.

One type of such constraints can be obtained as follows.
For a $0/1$-matrix $M$ with $n$ columns and a column subset $J \subseteq \{1, \ldots, n\}$, a \emph{minor} $M' = M \setminus U$ is obtained by 
\begin{itemize}[leftmargin=12pt, itemsep=0pt]
\item deleting a column $j \in J$ (that is, removing from $M$ column $j$ and all rows having a 1-entry in column $j$, which corresponds to setting $x_j = 1$ in the system $M \mathbf x \geq \mathbf 1$) or
\item contracting a column $j \in J$ (that is, removing from $M$ column $j$ and all rows $\mathbf y$ that became redundant as there is now another row $\mathbf y'$ with $\mathbf y' \leq \mathbf y$, which corresponds to setting $x_j = 0$ in the system $M \mathbf x \geq \mathbf 1$)
\end{itemize}
in any order. For such a minor $M' = M \setminus J$, the \emph{rank-constraint} associated with $M'$ is
\begin{equation}\label{eq_rank}
  x(M') =  \sum_{i \in \{1, \ldots, n\} \setminus J}x_i \geq \tau(M')
\end{equation}
and we have the following according to \cite{AB_2009}.

\begin{theorem}[\cite{AB_2009}] \label{thm_AB}
  Let $M$ be a $0/1$-matrix and $M'$ be some minor of $M$.
If the rank-constraint (\ref{eq_rank}) associated with $M'$ is a facet of $P(M')$, then it is also a facet of $P(M)$.
\end{theorem}

Hence, essential constraints associated with minors $M'$ remain essential for the covering polyhedron of the whole matrix $M$. 

For any hypergraph $\mathcal{H} = (V, \mathcal{E})$, let $M(\mathcal{H})$ denote its incidence matrix encoding row-wise its hyperedges $F$ (that is, the row of $M(\mathcal{H})$ corresponding to $F$ is a $0/1$-vector of length $|V|$ having a $1$-entry if $v \in F$ and a $0$-entry otherwise). Note that $\tau(\mathcal{H}) = \tau(M(\mathcal{H}))$ clearly holds.

The covering polyhedra associated with complete $q$-roses of order $n$ were studied in \cite{ABLW_2018,S_1989}.
In \cite{S_1989}, it has been proved that the rank-constraint associated with $M({\cal{R}}_n^q)$ defines a facet of $P(M({\cal{R}}_n^q))$.
In \cite{ABLW_2018}, it was shown that any minor of $M({\cal{R}}_n^q)$ corresponds to a complete $q'$-rose of order $n'$ with $q' \leq q$ and $n' \leq n$, and that no other constraints than rank-constraint associated with such minors define facets:

\begin{theorem}[\cite{ABLW_2018}] \label{thm_ABLW}
The covering polyhedron of the incidence matrix of a complete $q$-rose ${\cal{R}}_n^q$ of order $n$ is given by the non-negativity constraints and  
$$x(V') =  \sum_{v \in V'}x_v \geq |V'|-q +1$$
for all subsets $V' \subseteq \{1,\ldots,n\}$ with $|V'| \in \{q,\ldots,n\}$.
\end{theorem}
Note that this also reproves the result from \cite{S_1989} showing that $\tau({\cal{R}}_n^q)=n-q+1$ holds.

\subsection{Polyhedra associated with \OD-codes}
Following a similar approach as in e.g. \cite{ABLW_2018,ABLW_2022} for \X-problems with $\X \in \{\ID,\LD,\LTD,\OTD\}$, we study here polyhedra associated with \OD-codes. 
Due to $\od(G)=\tau(\hyp_{\OSD}(G)) = \tau(\clu_{\OSD}(G))$, we can determine a minimum \OD-code in a graph $G=(V,E)$ by solving the following covering problem
\[
\begin{array}{rcl}
\displaystyle \min \mathbf 1^T\mathbf x \ &   &    \\
M_{\OD}(G) \ \mathbf x \ & \geq & \mathbf 1   \\
\mathbf x \ & \in  & \mathbf Z_+^{|V|}  \\
\end{array}
\]
where $M_{\OD}(G)$ is the incidence matrix of the \OD-clutter $\clu_{\OD}(G)$. 
Accordingly, the \emph{\OD-polyhedron} of a graph $G$ is defined by 
$$
P_{\OD}(G)= P(M_{\OD}(G)) = \mbox{conv}\{\mathbf x \in \mathbf Z_+^{|V|} : M_{\OD}(G)\ \mathbf x \geq \mathbf 1\}. 
$$
Based on Theorem \ref{thm_BN}
from \cite{BN_1989} on general covering polyhedra, we prove the following.

\begin{theorem} \label{thm_OD-polyhedron}
Let $G=(V,E)$ be an \OD-admissible graph. $P_{\OD}(G)$ has
\begin{itemize}[leftmargin=18pt, itemsep=0pt]
\item[(a)] the equation $x_v = 1$ for all forced vertices $v \in \for_{\OD}(G)$;
\item[(b)] a non-negativity constraint $x_v \geq 0$ for all vertices $v \not\in \for_{\OD}(G)$; and
\item[(c)] a facet $\sum_{v \in F}x_v \geq 1$ for all hyperedges $F$ of $\clu_{\OD}(G)$ with $F \in \fac_{\OD}(G)$.
\end{itemize}
\end{theorem}

\begin{proof}
Let $G=(V,E)$ be an \OD-admissible graph. 
By the definition of $\for_{\OD}(G)$, we clearly have $x_v = 1$ for all forced vertices $v \in \for_{\OD}(G)$. 
This makes the non-negativity constraints for forced vertices redundant so that only $x_v \geq 0$ for all vertices $v \not\in \for_{\OD}(G)$ remain.

From Theorem \ref{thm_BN} by \cite{BN_1989} it is further known that the only facet-defining (i.e. essential) inequalities of a covering polyhedron $P(M)$ with integer coefficients and right hand side equal to $1$ are those of the defining system $M \mathbf x \geq \mathbf 1$. 
Hence, for $M = M_{\OD}(G)$, all constraints $\sum_{v \in F}x_v \geq 1$ for all hyperedges $F \in \fac_{\OD}(G)$ are facet-defining for $P_{\OD}(G)$.
\end{proof}

To illustrate these general properties, consider a half-graph $B_k=(U \cup W, E)$ for some $k \geq 2$ and recall from Lemma \ref{lem_OD-OTD_half-graphs} that $\for_{\OD}(B_k) = (U \setminus \{u_1\}) \cup (W \setminus \{w_k\})$ holds and that $\{u_1,w_k\}$ is the only hyperedge in $\fac_{\OD}(B_k)$. According to Theorem \ref{thm_OD-polyhedron}, $P_{\OD}(B_k)$ thus has 
\begin{itemize}[leftmargin=12pt, itemsep=0pt]
\item the equation $x_v = 1$ for all forced vertices $v \in \for_{\OD}(B_k) = (U \setminus \{u_1\}) \cup (W \setminus \{w_k\})$;
\item two non-negativity constraints $x_{u_1} \geq 0$ and $x_{w_k} \geq 0$; and 
\item one facet $x_{u_1} + x_{w_k}  \geq 1$ for the only hyperedge of $\clu_{\OD}(B_k)$ in $\fac_{\OD}(B_k)$.
\end{itemize}
This clearly shows that, on the one hand, $(U \setminus \{u_1\}) \cup W$ and $U \cup (W \setminus \{w_k\})$ are the only two minimum \OD-codes of $B_k$ and, on the other hand, that $P_{\OD}(B_k) = Q_{\OD}(B_k)$ holds and the incidence matrix of $\clu_{\OD}(B_k)$ is in fact ideal.

Determining the \OD-clutters of several graph families studied in Section \ref{sec4} showed their relation to complete $q$-roses ${\cal{R}}_n^q$ of order $n$. Relying on Theorem \ref{thm_ABLW} from \cite{ABLW_2018} enabled us to prove the following results. Recall that, for $q=2$, ${\cal{R}}_n^q$ is in fact the clique $K_n$.

\begin{theorem} \label{thm_OD-polyhedron_Kn}
Let $G=(V,E)$ be either a clique $K_n$ with $n \geq 2$ or a matching $kK_2$ with $k \geq 1$ and $n=2k$. 
Then, we have $\clu_{\OD}(G)= {\cal{R}}_n^2=K_n$ and $P_{\OD}(G)$ is given by
\begin{itemize}[leftmargin=18pt, itemsep=0pt]
\item[(a)] a non-negativity constraint $x_v \geq 0$ for all vertices $v \in V$; and
\item[(b)] $x(V') =  \sum_{v \in V'}x_v \geq |V'|-1$ for all subsets $V' \subseteq V$ with $|V'| \geq 2$.
\end{itemize}
\end{theorem}

\begin{proof}
Let $G=(V,E)$ be either a clique $K_n$ with $n \geq 2$ or a matching $kK_2$ with $k \geq 1$ and $n=2k$. 
From Lemma \ref{lem_OD_cliques} and the proof of Lemma \ref{lem_OD_union-cliques}, respectively, 
we see that $\clu_{\OD}(K_n)= {\cal{R}}_n^2=K_n$ and $\clu_{\OD}(kK_2)= {\cal{R}}_{2k}^2=K_{2k}$ holds. 
Hence, $\clu_{\OD}(G)$ is in both cases a complete $2$-rose of order $n$ and $P_{\OD}(G)$ is accordingly given by 
non-negativity constraints for all vertices $v \in V$ and 
constraints $x(V') \geq |V'|-1$ for all subsets $V' \subseteq V$ with $|V'| \geq 2$ by Theorem \ref{thm_ABLW} from \cite{ABLW_2018}.
\end{proof}

Note that two graphs with equal \OD-clutters have the same set of \OD-codes and thus also the same \OD-numbers and \OD-polyhedra. Theorem~\ref{thm_OD-polyhedron_Kn} shows that this applies to cliques and matchings.

Furthermore, consider a $k$-fan $F_k = (W_2 \cup \{u\},E)$ for some $k \geq 2$ and recall from the remark after the proof of Theorem \ref{thm_OD_clique-stars} that $\clu_{\OD}(F_k)=K_{2k}$ (where the $K_{2k}$ is induced by $W_2$ and the universal vertex $u$ of $F_k$ becomes an isolated vertex in $\clu_{\OD}(F_k)$ not contained in any hyperedge). Similar arguments as in the proof of Theorem \ref{thm_OD-polyhedron_Kn} show that $P_{\OD}(F_k)$ has non-negativity constraints for all vertices, and constraints $x(V') \geq |V'|-1$ for all subsets $V' \subseteq W_2$ with $|V'| \geq 2$.

The following two theorems show that the \OD-numbers of thin and thick headless spiders, as calculated in Lemma \ref{lem_ThickSpider_OD_clutter} and \ref{lem_ThinSpider_OD_clutter}, 
can also be deduced by the use of polyhedral techniques.

\begin{theorem} \label{thm_OD-polyhedron_coHk}
Let $\overline H_k=(Q \cup S, E)$ be a thick headless spider with $k \geq 4$. Then, we have $\clu_{\OD}(\overline H_k) = {\cal{R}}_{|S|}^{|S|-1} \cup {\cal{R}}_{|Q|}^2$ and $P_{\OD}(\overline H_k)$ is given by the constraints 
\begin{itemize}[leftmargin=18pt, itemsep=0pt]
\item[(a)] $x_v \geq 0$ for all vertices $v \in Q \cup S$,
\item[(b)] $x(V') =  \sum_{v \in V'}x_v \geq |V'|-k+2$ for all $V' \subseteq S$ with $|V'| \geq k-1$,
\item[(c)] $x(V') =  \sum_{v \in V'}x_v \geq |V'|-1$ for all $V' \subseteq Q$ with $|V'| \geq 2$.
\end{itemize}
\end{theorem}

\begin{proof}
Consider a thick headless spider $\overline H_k=(Q \cup S, E)$ with $k \geq 4$. 
From Lemma \ref{lem_ThickSpider_OD_clutter}, we have $\clu_{OD}(\overline H_k) = {\cal{R}}_{|S|}^{|S|-1} \cup {\cal{R}}_{|Q|}^2$.
By $\for_{\OD}(\overline H_k) = \emptyset$, we have non-negativity constraints for all vertices. 
Applying Theorem \ref{thm_ABLW} from \cite{ABLW_2018} on polyhedra associated to complete $q$-roses shows that $P_{\OD}(\overline H_k)$ is given by non-negativity constraints for all vertices and the constraints 
\begin{itemize}[leftmargin=12pt, itemsep=0pt]
\item $x(V') =  \sum_{v \in V'}x_v \geq |V'|-k+2$ for all $V' \subseteq S$ with $|V'| \geq k-1$,
\item $x(V') =  \sum_{v \in V'}x_v \geq |V'|-1$ for all $V' \subseteq Q$ with $|V'| \geq 2$.
\end{itemize}
\end{proof}
Note that the facets $x(S) \geq 2$ and $x(Q) \geq k-1$ together also imply that $\od(\overline H_k) = k+1$ holds. 

Comparing Lemma \ref{lem_ThickSpider_OD_clutter} and Theorem \ref{thm_OD-polyhedron_coHk}
with the result from \cite{ABLW_2022} on \OTD-codes of thick headless spiders, we observe that $\clu_{\OD}(\overline H_k) = \clu_{\OTD}(\overline H_k)$. Hence, a vertex subset is an \OD-code of $\overline H_k$ if and only if it is an \OTD-code of $\overline H_k$. 
Accordingly, the \OD- and \OTD-numbers as well as the \OD- and \OTD-polyhedra are equal for thick headless spiders.

\begin{theorem} \label{thm_OD-polyhedron_Hk}
Consider a thin headless spider $H_k=(Q \cup S, E)$ with $k \geq 4$. Then, we have $\clu_{\OD}(H_k) = H_k$ and $P_{\OD}(H_k)$ is given by the constraints 
\begin{itemize}[leftmargin=18pt, itemsep=0pt]
\item[(a)] $x_v \geq 0$ for all vertices $v \in Q \cup S$,
\item[(b)] $x_{q_i} + x_{s_i} \geq 1$ for $1 \leq i \leq k$,
\item[(c)] $x(V') =  \sum_{v \in V'}x_v \geq |V'|-1$ for all $V' \subseteq Q$ with $|V'| \geq 2$.
\end{itemize}
\end{theorem}

\begin{proof}
Consider a thin headless spider $H_k=(Q \cup S, E)$ with $k \geq 4$.
From Lemma \ref{lem_ThinSpider_OD_clutter}, we have $\clu_{\OD}(H_k) = H_k$. 
Hence $\mathcal{F}_{\OD}(H_k)$ is composed of the matching $\{q_i,s_i\}$ for $1 \leq i \leq k$ and the complete $2$-rose ${\cal{R}}_{|Q|}^2$. 

By $\for_{\OD}(H_k) = \emptyset$, we have non-negativity constraints for all vertices. 
Moreover, $P_{\OD}(H_k)$ has clearly the constraints $x(V') \geq |V'|-1$ for all $V' \subseteq Q$ with $|V'| \geq 2$ by Theorem \ref{thm_ABLW} from \cite{ABLW_2018} and Theorem \ref{thm_AB} from \cite{AB_2009}. 
Finally, $x_{q_i} + x_{s_i} \geq 1$ define facets for $1 \leq i \leq k$ by Theorem \ref{thm_BN} from \cite{BN_1989}, 
and it is easy to see that no further constraints are needed to describe $P_{\OD}(H_k)$ (as each $x_{s_i}$ occurs in exactly one constraint different from a non-negativity constraint).
\end{proof}

Note that combining all the constraints in Theorem \ref{thm_OD-polyhedron_Hk}(b) yields $x(Q) + x(S) \geq k$ and this implies $\gamma^{\OSD}(H_k) \geq k$. 
It is also easy to see that $Q$ is a cover of $\clu_{\OD}(H_k)$ and hence, $\od(H_k) = k$.

Furthermore, consider an extended thin spider $E_k=(Q \cup S, E)$ with $k \geq 4$. 
Lemma \ref{lem_ExtThinSpider_OD-OTD_clutter} shows that $\clu_{\OD}(E_k)$ differs from $\clu_{\OD}(H_k)$ in exactly one hyperedge, namely $\{s_k\}$ replaces $\{q_k,s_k\}$, whereas all other hyperedges $\{q_i,s_i\}$ for all $1 \leq i < k$, and $\{q_i,q_j\}$ for all $1 \leq i < j \leq k$ are the same in both cases. Accordingly, $P_{\OD}(E_k)$ can be obtained from $P_{\OD}(H_k)$ by replacing the constraint $x_{q_k} + x_{s_k} \geq 1$ by the equation $x_{s_k} = 1$.

Finally, recall that thin headless spiders can be further generalized to thin suns $T_k=(C\cup S,E)$ with $k \geq 3$.
We next study the \OD-polyhedra of two families of thin suns, namely 
sunlets and almost complete thin suns as examples of thin suns without and with open $C$-twins, respectively.

We start with sunlets and recall that $T_4$ is the only sunlet with open $C$-twins, hence we study the case with $k \geq 5$.

\begin{lemma}\label{lem_P-Sunlet} 
  For a sunlet $T_k=(C \cup S, E)$ with $k \geq 5$, $P_{\OD}(T_k)$ has non-negativity constraints for all vertices and the constraints
\begin{itemize}[leftmargin=12pt, itemsep=0pt]
\item $x(V') =  \sum_{v \in V'}x_v \geq |V'|-1$ for all $V' \subseteq \{s_i,c_{i-1},c_i,c_{i-+}\}$ with $|V'| \geq 2$ and $1 \leq i \leq k$,
\item $x(C') =  \sum_{v \in C'}x_v \geq |C'|-1$ for all $C' \subseteq C$ with $|C'| \geq 2$.
\end{itemize}
\end{lemma}

\begin{proof}
  Consider a sunlet $T_k=(C \cup S, E)$ with $k \geq 5$.
  By Lemma \ref{lem_ThinSun_clutter}, the \OD-clutter $\clu_{\OD}(T_k)$ is composed of
\begin{itemize}[leftmargin=12pt, itemsep=0pt]
\item $N[s_i] = \{s_i,c_i\}$ for all $s_i \in S$,
\item $\Delta(c_i,s_{i-1}) = \{s_i,c_{i-1}\}$ and $\Delta(c_i,s_{i+1}) = \{s_i,c_{i+1}\}$ for all $c_i \in C$ (since we have $\Delta(c_i,s_{i \pm 1}) = \{s_i,c_{i-1},\\ c_{i+1}\} \Delta \{c_{i \pm 1}\}$), 
\item $\Delta(s_i,s_j) = \{c_i,c_j\}$ for all distinct $s_i,s_j \in S$. 
\end{itemize}
Hence, in $\clu_{\OD}(T_k)$, $C$ induces a $K_k$, and all vertex subsets $\{s_i,c_{i-1},c_i,c_{i+1}\}$ induce a $K_4$ for $1 \leq i \leq k$.
Accordingly, $\for_{\OD}(T_k) = \emptyset$ implies that $P_{\OD}(T_k)$ has non-negativity constraints for all vertices.
Moreover, we deduce from Theorem \ref{thm_ABLW} by \cite{ABLW_2018} and Theorem \ref{thm_AB} from \cite{AB_2009} that $P_{\OD}(T_k)$ has the constraints
\begin{itemize}[leftmargin=12pt, itemsep=0pt]
\item $x(V') =  \sum_{v \in V'}x_v \geq |V'|-1$ for all $V' \subseteq \{s_i,c_{i-1},c_i,c_{i-+}\}$ with $|V'| \geq 2$ and $1 \leq i \leq k$,
\item $x(C') =  \sum_{v \in C'}x_v \geq |C'|-1$ for all $C' \subseteq C$ with $|C'| \geq 2$.
\end{itemize}
\end{proof}

We finally consider almost complete thin suns and can show the following.

\begin{lemma}\label{lem_P-almost-complete-thin-sun}
For an almost complete thin sun $T_{2\ell}$ with $\ell \geq 3$, $P_{\OD}(T_{2\ell})$ has non-negativity constraints for all vertices and the constraints
\begin{itemize}[leftmargin=12pt, itemsep=0pt]
\item $x_{s_i} + x_{s_{i+\ell}} \geq 1$ for $1 \leq i \leq \ell$,
\item $x_{s_i} + x_{c_i} \geq 1$ for $1 \leq i \leq 2\ell$,
\item $x(C') =  \sum_{v \in C'}x_v \geq |C'|-1$ for all $C' \subseteq C$ with $|C'| \geq 2$.
\end{itemize}
\end{lemma}

\begin{proof}
Consider an almost complete thin sun $T_{2\ell}$ with $\ell \geq 3$ where $c_i,c_{i+\ell}$ form open $C$-twins. 
Recall that, by Lemma \ref{lem_ThinSun_clutter}, the \OD-clutter $\clu_{\OD}(T_{2\ell})$ is composed of 
\begin{itemize}[leftmargin=12pt, itemsep=0pt]
\item $N[s_i] = \{s_i,c_i\}$ for all $s_i \in S$,
\item $\Delta(s_i,s_j) = \{c_i,c_j\}$ for all distinct $s_i,s_j \in S$, 
\item $\Delta(c_i,c_{i+\ell}) = \{s_i,s_{i+\ell}\}$ for all $c_i \in C$.
\end{itemize}
Hence, in $\clu_{\OD}(T_{2\ell})$, $C$ induces a $K_k$, and all vertex subsets $\{s_i,c_i,c_{i+ \ell},s_{i+\ell}\}$ induce a $C_4$ for $1 \leq i \leq \ell$. 
Accordingly, $\for_{\OD}(T_{2\ell}) = \emptyset$ implies that $P_{\OD}(T_{2\ell})$ has non-negativity constraints for all vertices.
Moreover, $P_{\OD}(T_{2\ell})$ has the constraints
\begin{itemize}[leftmargin=12pt, itemsep=0pt]
\item $x_{s_i} + x_{s_{i+\ell}} \geq 1$ for $1 \leq i \leq \ell$,
\item $x_{s_i} + x_{c_i} \geq 1$ for $1 \leq i \leq 2\ell$,
\end{itemize}
from the defining system and by Theorem \ref{thm_ABLW} from \cite{ABLW_2018} and Theorem \ref{thm_AB} from \cite{AB_2009} in addition the constraints
\begin{itemize}[leftmargin=12pt, itemsep=0pt]
\item $x(C') =  \sum_{v \in C'}x_v \geq |C'|-1$ for all $C' \subseteq C$ with $|C'| \geq 2$.
\end{itemize}
\end{proof}

We conjecture in both cases for sunlets and almost complete thin suns that the \OD-polyhedra have no other facets than those presented in the two previous lemmas.

The results presented in this section illustrates how, on the one hand, polyhedral arguments can be used to determine lower bounds for \OD-numbers and, on the other, an analysis of the \OD-clutter provides \OD-codes. Moreover, if the order of the latter meets the lower bound, the \OD-number of the studied graph is determined.

We note further that manifold hypergraphs and matrices have been already studied in the covering context, see for instance 
\cite{AB_2009,C_2001,NS_1989,S_1989} to mention just a few. 
The same techniques as illustrated above with the help of complete $q$-roses can be applied whenever the \OD-clutter of some graph equals such a hypergraph or contains it as a minor, which gives an interesting perspective of studying \OD-polyhedra further.

\section{Concluding remarks}

In this paper, we introduced and studied open-separating dominating codes in graphs. 
We showed that such codes exist in graphs without open twins and that finding minimum \OSD-codes is NP-hard. This also enables us to reprove the \NP-hardness of the \OTD-problem on a narrower class of input graphs (bipartite, four-cycle-free and of girth at least~$6$) and with an arguably simpler reduction than previously considered, for example, in~\cite{PP_2017,SS2010}. Moreover, we provided bounds on the \OSD-number of a graph both in terms of its number of vertices and in relation to other \X-numbers, notably showing that \OSD- and \OLD-number of a graph differ by at most one. Despite this closeness between the \OSD- and the \OLD-numbers, we proved that it is \NP-hard to decide if the two said parameters of a graph actually differ. This further motivated us to compare the two numbers on several graph families. This study revealed that they 
\begin{itemize}
\item are equal, for example, for cliques $K_n$ with $n \geq 3$, 
thin and thick headless spiders $H_k$ and $\overline H_k$, respectively, with $k \geq 3$, and 
thin suns $T_k=(C \cup S, E)$ with $k \geq 4$ and without open $C$-twins;
\item differ for example, for matchings $kK_2$ with $k \geq 1$, $k$-fans $F_k$ with $k \geq 2$, 
half-graphs $B_k$ with $k \geq 1$ and their disjoint unions, $k$-double stars $D_k$ with $n \geq 2$, extended thin spiders $E_k$ with $k \geq 4$ and 
almost complete thin suns $T_{2\ell}$ with $\ell \geq 3$.
\end{itemize}
In particular, this showed that the \OSD-numbers of cliques, half-graphs and their disjoint unions attain the upper bound in Theorem \ref{thm_osd_bounds-n}. In addition, matchings and disjoint unions of half-graphs are examples of disconnected graphs $G$ where $\od(G)$ is larger than the sum of the \OD-numbers of its components. This behavior does not apply to any other \X-number for $\X \in \{\ID, \DTD, \LD, \LTD, \OLD\}$.

Moreover, we provided an equivalent reformulation of the \OSD-problem as a covering problem in a suitable hypergraph composed of the closed neighborhoods and the symmetric differences of open neighborhoods of vertices. We also discussed the polyhedra associated with the \OSD-codes, particularly, in relation to the \OLD-codes of some graph families already studied in this context. 
The latter illustrated how polyhedral arguments can be used to determine lower bounds for \OSD-numbers, how an analysis of the \OSD-clutter can provide the \OSD-codes, and that combining both arguments can yield the \OSD-numbers of the studied graphs.

The future lines of our research include studying the \OSD-problem on more graph families and also searching for extremal cases concerning the lower bounds for \OSD-numbers (that is, the logarithmic bound in Theorem \ref{thm_osd_bounds-n} in terms of the order of the graph and the lower bounds in Theorem \ref{thm_osd_bounds-L} based on \LD- and \LTD-numbers). Even though the problem of deciding if the \OSD- and the \OLD-numbers differ is \NP-complete in general, it would be interesting to see if for some particular graph families, this problem becomes polynomial-time solvable. In that case, it would be further interesting to provide a complete dichotomy as to for which graphs of that graph family the two code-numbers differ and for which they are equal.

\section*{Acknowledgement}
This research was financed by a public grant overseen by the French National Research Agency as part of the ``Investissements d’Avenir'' through the IMobS3 Laboratory of Excellence (ANR-10-LABX-0016), by the French government IDEX-ISITE initiative 16-IDEX-0001 (CAP 20-25) and the International Research Center ``Innovation Transportation and Production Systems'' of the I-SITE CAP 20-25.
\bibliographystyle{splncs04}

\end{document}